\RequirePackage{fix-cm}
\pdfoutput=1
\documentclass[leqno, fleqn, centertags, 12pt]{article}

\usepackage{latexsym}
\usepackage{amsmath}
\usepackage{amscd}
\usepackage{amssymb}
\usepackage{verbatim}

\usepackage[T1]{fontenc}
\usepackage{latexsym}
\usepackage{manfnt}
\usepackage[sans]{dsfont}
\usepackage{palatino}
\usepackage[sc, osf]{mathpazo} 
\usepackage{eulervm}
\usepackage{eufrak}
\usepackage{euscript}

\usepackage{textcomp}

\usepackage{fullpage}

\usepackage{graphicx}

\newcommand{\id}{\mathop{{\mbox{id}}\fff}\nolimits}

\newcommand{\gen}{{\mathfrak g}}
\newcommand{\genus}{\mathop{\mathfrak g}\hspace{0.075em}}
\newcommand{\mmod}{\mathop{\mbox{Mod}}\hspace{0.075em}}
\newcommand{\mm}{\mathop{\mbox{Mod}}\hspace{0.075em}}
\newcommand{\mms}{\mmod(S)}
\newcommand{\ttt}{\mathop{\EuScript I}}
\newcommand{\tors}{\ttt(S)}

\newcommand{\rank}{\mathop{\mbox{rank}}\hspace{0.1em}}
\newcommand{\rankf}{\mathop{\mbox{\emph{rank}}}\hspace{0.1em}}

\newcommand{\zzz}{\mathbold{Z}}

\newcommand{\toto}{\longrightarrow}
\newcommand{\ttoo}{\hspace*{0.2em}\longrightarrow\hspace*{0.2em}}

\def\sss{\hspace{0.05em}\ }
\def\dss{\hspace{0.1em}\ }
\def\trs{\hspace{0.15em}\ }
\def\qss{\hspace{0.2em}\ }
\def\pss{\hspace{0.3em}\ }
\def\oss{\hspace{0.4em}\ }

\def\halfff{\hspace*{0.025em}}
\def\fff{\hspace*{0.05em}}
\def\dff{\hspace*{0.1em}}
\def\trf{\hspace*{0.15em}}
\def\qff{\hspace*{0.2em}}
\def\pff{\hspace*{0.3em}}
\def\off{\hspace*{0.4em}}

\def\ttff{{\hspace*{0.1em}--\hspace*{0.1em}}}

\def\dfdot{\hspace*{-0.2em}.\hspace*{0.4em}\ }

\def\dfcom{\hspace*{-0.2em},\hspace*{0.4em}\ }

\newcommand{\nsp}{\hspace*{-0.1em}}

\newcommand{\snsp}{\hspace*{-0.175em}}
\newcommand{\dnsp}{\hspace*{-0.2em}}

\newcommand{\eeq}{ \qff =\qff }

\newcommand{\ccomp}{\mathsf{c}}

\newcommand{\ccut}{\nsp/\!\!/\hspace{-0.08em}}
\newcommand{\cutc}{\nsp/\!\!/\hspace{-0.08em}c}

\newcommand{\hclass}[1]{[\dff #1\dff]}

\makeatletter
\renewcommand{\@makefntext}[1]{\vspace*{0.5ex}\parindent=0em
\hspace*{-0.4em}
\hbox to 0.4em{\hss\@makefnmark}\hspace*{0.4em}{#1}
}
\makeatother

\newcommand{\proof}{\vspace{0.5\bigskipamount}{\textbf{{\emph{Proof}.}}\hspace*{0.7em}}}
\newcommand{\prooftitle}[1]{\vspace*{0.5\bigskipamount}{\textbf{{\emph{#1}.}}\hspace{0.7em}}}
\newcommand{\eproof}{ $\blacksquare$}
\newcommand{\esubproof}{ $\square$}

\newcounter{mysectionnumber}
\newcommand{\mysection}[2]{
\setcounter{equation}{0}\refstepcounter{mysectionnumber}
\section*{ \textnormal{{\themysectionnumber.} {#1}}}\label{#2}}

\newcounter{myparnum}[mysectionnumber]
\newcommand{\mypar}[2]{\refstepcounter{myparnum}{\vspace{\medskipamount}
\textbf{\textit{{\themyparnum. #1}\label{#2}}}\hspace{0.5em}}}
\renewcommand{\themyparnum}{\themysectionnumber.\arabic{myparnum}}

\numberwithin{equation}{section}

\newcommand{\myitpar}[1]{\vspace{\medskipamount}\textbf{\textit{#1}}\hspace*{0.5em}}
\newcommand{\myit}[1]{\textbf{\textit{#1}}\hspace{0.0em}}

\newcommand{\mynonumbersection}[2]{
\vspace{-0.0ex}
\section*{{}\hspace*{0.00em}$\phantom{1.}$ \textnormal{{#1}}}\label{#2}\vspace{\bigskipamount}\vspace{-0.0ex}}

\newcounter{myappendnumber}
\newcommand{\myappend}[2]{\setcounter{footnote}{0}
\setcounter{myapparnum}{0}
\refstepcounter{myappendnumber}
\section*{\textnormal{Appendix.\oss 
{#1}}}\label{#2}}

\newcounter{myapparnum}[mysectionnumber]
\newcommand{\myappar}[2]{\refstepcounter{myapparnum}{\vspace{\bigskipamount}
\textbf{\textit{{\themyapparnum. #1}\label{#2}}}\hspace{0.5em}}}
\renewcommand{\themyapparnum}{A.
\arabic{myapparnum}}

\title{Elements\qss of\qss Torelli\qss topology\halfff:\oss \\ 
I.\qss The\qss rank\qss of\qss abelian\qss subgroups}
\author{Nikolai\dss V.\dss Ivanov}
\date{}

\begin{document}

\setlength{\baselineskip}{12pt plus 0pt minus 0pt}
\setlength{\parskip}{12pt plus 0pt minus 0pt}
\setlength{\abovedisplayskip}{12pt plus 0pt minus 0pt}
\setlength{\belowdisplayskip}{12pt plus 0pt minus 0pt}

\newskip\smallskipamount \smallskipamount=3pt plus 0pt minus 0pt
\newskip\medskipamount   \medskipamount  =6pt plus 0pt minus 0pt
\newskip\bigskipamount   \bigskipamount =12pt plus 0pt minus 0pt

\maketitle

\vspace*{12ex}

\myit{\hspace*{0em}\large Contents}  \vspace*{1ex} \vspace*{\bigskipamount}\\ 
\myit{Preface}\hspace*{0.5em}  \hspace*{0.5em} \vspace*{1ex}\\
\myit{\phantom{1}1.}\hspace*{0.5em} Surfaces, circles, and diffeomorphisms  \hspace*{0.5em} \vspace*{0.25ex}\\
\myit{\phantom{1}2.}\hspace*{0.5em} Dehn twists  \hspace*{0.5em} \vspace*{0.25ex}\\
\myit{\phantom{1}3.}\hspace*{0.5em} Action of Dehn twists on homology  \hspace*{0.5em} \vspace*{0.25ex}\\
\myit{\phantom{1}4.}\hspace*{0.5em} Dehn multi-twists in Torelli groups  \hspace*{0.5em} \vspace*{0.25ex}\\
\myit{\phantom{1}5.}\hspace*{0.5em} The rank of multi-twist subgroups of Torelli groups \hspace*{0.5em} \vspace*{0.25ex}\\  
\myit{\phantom{1}6.}\hspace*{0.5em} Pure diffeomorphisms and reduction systems  \hspace*{0.5em} \vspace*{0.25ex}\\ 
\myit{\phantom{1}7.}\hspace*{0.5em} The rank of abelian subgroups of Torelli groups  \hspace*{0.25em} \vspace*{0.25ex}\\
\myit{\phantom{1}8.}\hspace*{0.5em} Dehn and Dehn--Johnson twists in Torelli groups\halfff:\qss I  \hspace*{0.5em} \vspace*{0.25ex}\\
\myit{\phantom{1}9.}\hspace*{0.5em} Dehn and Dehn--Johnson twists in Torelli groups\halfff:\qss II  \hspace*{0.5em} \vspace*{1ex}\\
\myit{Appendix.}\hspace*{0.5em} Necklaces  \hspace*{0.5em} \vspace*{1ex}\\
\myit{References}\hspace*{0.5em}  \hspace*{0.5em}  \vspace*{0.25ex}

\footnotetext{\hspace*{-0.65em}\copyright\qss 
Nikolai\qss V.\qss Ivanov,\oss 2016.\oss 
Neither the work reported in this paper,\qss 
nor its preparation were 
supported by any governmental 
or non-governmental agency,\qss 
foundation,\qss 
or institution.}

\renewcommand{\baselinestretch}{1}
\selectfont

\newpage
\mynonumbersection{Preface}{preface}

\vspace*{\bigskipamount}
By Torelli topology the author understands aspects of the topology of surfaces\qss 
(potentially)\qss relevant to the study of Torelli groups.\oss

The present paper is devoted to a new approach to the results of\qss W.\qss Vautaw\qss \cite{v}\qss
about Dehn multi-twists in Torelli groups and abelian subgroups of Torelli groups.\oss
The main results are a complete description of Dehn multi-twists in Torelli groups
and the theorem to the effect that the rank of any abelian subgroup of the Torelli group
of a closed surface of genus\qss $\gen$\qss is\qss $\leq\qff 2\dff\gen\qff -\qff 3$\snsp.\oss
In contrast with W.\qss Vautaw's paper\qss \cite{v},\oss 
which heavily relies on the graph-theoretic language,\oss
the present paper is based on topological methods.\oss
The resulting proofs are more transparent and lead to stronger estimates
of the rank of an abelian subgroup when some additional information is available.\oss
A key role is played by the notion of a\qss 
\emph{necklace}\trs of a system of circles on a surface.\oss

As an unexpected application of our methods,\qss 
in Sections\dss \ref{simplest-twists-sufficient}\dss 
and\dss \ref{simplest-twists-necessary}\dss we give a new proof of
the algebraic characterization of the Dehn twist about separating circles
and the Dehn-Johnson twists about bounding pairs of circles from the paper\qss \cite{fi}\qss
of\qss B.\qss Farb and the author\halfff.\oss
This proof is much shorter than the original one and bypasses one of the main
difficulties,\oss which may be called the\qss \emph{extension problem},\qss 
specific to Torelli groups as opposed 
to the Teichm\"{u}ller modular groups.\oss
This algebraic characterization is also contained
in the Ph.D.\qss Thesis of W.\qss Vautaw\qss \cite{vt}.\oss
But the proof in\qss \cite{vt}\qss is not quite complete because the extension problem 
is ignored.\oss
The extension problem is discussed in details in the second paper\qss 
\cite{i-bound}\qss of this series.\oss

In many respects Torelli groups look similar to
Teichm{\"u}ller modular groups from a sufficient distance.\oss
For example,\oss the algebraic characterization of the Dehn twist about separating circles
and the Dehn-Johnson twists about bounding pairs of circles in Torelli groups
is morally the same as author's\qss \cite{i-aut}\qss
algebraic characterization of Dehn twists in Teichm\"{u}ller modular groups.\oss
But the analogy breaks down to a big extent if one focuses on the details.\oss
The extension problem is one of such details.\oss

The new proof this algebraic characterization is simpler thanks
not to any technical advances,\oss but to a few shifts in the point of view.\oss
It may well be the case that key new insight is 
in Theorem\qss \ref{multi-twists-subgroups}\qss below,\oss 
which is a trivial corollary of our form of the description\qss
of Dehn multi-twists in the Torelli groups\qss
(see\dss Theorem\qss \ref{multi-twists-in-t}\fff).\oss
Another novel aspect of this proof is the shift of emphasis
from the canonical reduction systems to the\qss
pure reduction systems.\qss
See Section\qss \ref{pure-reduction}.\oss

For the experts,\oss one perhaps should point out that the centers of the centralizers,\pss
a favorite tool of the author\halfff,\pss
are replaced by more elegant\qss \emph{bicommutants}.\oss
On the technical level these notions are equivalent 
in the situation of the present paper\halfff,\oss 
but bicommutants lead to a slightly different mindset\halfff.\oss 

Sections\qss \ref{prelim}{\ttff}\ref{dehn-torelli}\qss are devoted mostly
to a review of terminology and prerequisites,\pss
except of an instructive example in Section\qss \ref{dehn-torelli}.\oss
To a big extent this is also true for Section\qss \ref{pure-reduction},\oss
devoted to the reduction systems,\pss
except of rephrasing some well known results 
in the terms of pure reduction systems.\pss
The Appendix is devoted to some motivation for the term\trs \emph{necklace},\pss
and is not used in the main part of the paper\halfff.\pss
The table of contents should serve as a sufficient guide for the rest of the paper\halfff.\oss

\mysection{Surfaces,\qss circles,\qss and\qss diffeomorphisms}{prelim}

\myitpar{Surfaces.} By a\qss \emph{surface}\qss we understand a compact orientable $2$-manifold 
with\qss (possibly empty)\qss boundary.\oss
The boundary of a surface\qss $S$\qss is denoted by\qss $\partial S$\dnsp,\oss 
and its genus by\qss $\genus(S)$\dnsp.\oss
If\qss $\partial S\eeq\varnothing$\dfcom then\qss $S$\qss is called a\qss \emph{closed}\qss surface.\oss
By a\qss \emph{subsurface}\qss of a surface $S$ we understand a codimension $0$ submanifold\qss
$Q$\qss of\qss $S$\qss such that each component of\qss $\partial Q$\qss is either equal to
a component of\qss $\partial S$\dfcom or disjoint from\qss $\partial S$\dfdot 

For a subsurface $Q$ of $S$ we will denote by $\ccomp Q$ 
the closure of its set-theoretic complement $S\smallsetminus\nsp Q$\dfdot
Clearly,\oss $\ccomp Q$\qss is also a subsurface of\qss $S$\halfff\dfcom 
and\qss $\partial Q\smallsetminus \partial S\eeq\partial \ccomp Q\smallsetminus \partial S$\dfdot
We will say that\qss $\ccomp Q$\qss is the subsurface\dss \emph{complementary to}\dss $Q$.

\myitpar{Teichm\"{u}ller modular groups and Torelli groups.} The\qss 
\emph{Teichm\"{u}ller modular group}\qss 
$\mms$\qss of an orientable surface\qss $S$\qss is the group of isotopy classes of
orientation-preserving diffeomorphisms\qss $S\ttoo S$\dnsp.\oss
Both diffeomorphisms and isotopies are required to preserve the boundary\qss $\partial S$\qss
only set-wise\qss (this is automatic for diffeomorphisms and standard for isotopies).\oss
For closed surfaces\dss $S$\dss the subgroup of elements\qss $\mmod(S)$\qss
acting trivially on the homology group\qss $H_1(S\fff,\pff \zzz)$\qss
is called the\qss \emph{Torelli\dss group}\qss of\dss $S$\dss
and is denoted by\qss $\tors$\dnsp.\oss
Usually we will denote the homology group\qss $H_{\fff 1}(S\fff,\pff \zzz)$\qss
simply by\qss $H_{\fff 1}(S)$\dnsp.\oss

There are several candidates for the definition of 
the Torelli groups of surfaces with boundary\nsp,\oss
but none of them is completely satisfactory\nsp.\oss
It seems that this is so not for the lack of trying 
to find a\qss ``right definition'',\oss
but because there is no such definition.\oss

By technical reasons we will need also the subgroups\pss $\ttt_{\fff m\fff}(S)$\pss of\qss $\mms$\qss
consisting of the isotopy classes of diffeomorphisms acting trivially on\qss
$H_{\fff 1}(S\fff,\pff \zzz/m\dff\zzz)$\dnsp,\oss
where\qss $m$\qss is an integer\halfff.\oss
Clearly,\oss $\tors\qff\subset\qff\ttt_{\fff m\fff}(S)$\pss for every\qss $m\dff\in\dff \zzz$\snsp.\oss

\myitpar{Circles.} A\dss \emph{circle}\qss on a surface\dss $S$\dss is defined as
a submanifold of\dss $S$\dss diffeomorphic to 
the standard circle\dss $S^1$\dss and disjoint from\dss $\partial S$\dfdot
A circle in\dss $S$\dss is called\qss \emph{non-peripheral}\qss if does not bound 
an annulus together with a component of the boundary\dss $\partial S$\dfdot 
If\dss $S$\dss is a closed surface,\oss then,\oss obviously,\oss all circles in $S$ are non-peripheral.\oss 
A circle on\dss $S$\dss is said to be\qss \emph{non-trivial}\qss 
in\dss $S$\dss if it is if does not bound 
an annulus together with a component of the boundary\dss $\partial S$\dss
and does not bound a disc in\dss $S$\dfdot

\myitpar{Separating circles.} A circle $D$ in a connected surface $S$ 
is called\qss \emph{separating}\qss if the set-theoretic difference\dss 
$S\smallsetminus D$\dss is not connected.\oss 
In this case\dss $S\smallsetminus D$\dss consist of two components.\oss 
The closures of these components are subsurfaces of $S$ having $D$ as a boundary component.\oss
The other boundary components of these subsurfaces 
are at the same time boundary components of $S$\dfdot
We will call these subsurfaces\qss 
\emph{the parts into which\dss $D$\dss divides}\qss $S$\dfdot
If $Q$ is one of these parts,\oss then,\oss obviously,\oss $\ccomp Q$\qss is the other part.

\myitpar{Bounding pairs of circles.} A\qss \emph{bounding pair of circles}\qss 
on a connected surface $S$ is defined as an unordered pair\dss 
$C\fff,\pff C'$\dss of disjoint non-isotopic circles in $S$
such that both circles\dss $C\fff,\pff C'$\dss are non-separating\halfff,\oss 
but $S\smallsetminus (C\cup C')$ is not connected.\oss
In this case the difference $S\smallsetminus (C\cup C')$ consist of two components.\qss
The closures of these components are subsurfaces of 
$S$ having both $C$ and $C'$ as their boundary components.\oss
The other boundary components of these subsurfaces 
are at the same time boundary components of $S$\dfdot 
We will call these subsurfaces\qss 
\emph{the parts into which the bounding pair of circles\dss $C\fff,\pff C'$\dss divides}\qss $S$\dfdot
If\qss $Q$\qss is one of these parts,\oss 
then,\oss obviously,\oss $\ccomp Q$ is the other part.\oss

\myitpar{Cutting surfaces and diffeomorphisms.}
Every one-dimensional closed submanifold $c$ of a surface $S$ leads to a 
new surface\qss $S\cutc$\qss
obtained by cutting\qss $S$\qss along\qss $c$\dfdot 
The components of\qss $S\cutc$\qss are called the\qss \emph{parts}\qss
into which $c$ divides\qss $S$\dnsp.\oss
The canonical map\oss 
\[
\quad
p\cutc\dff\colon\dff S\cutc\off\toto\off S\fff.
\] 
induces a diffeomorphism\oss
$\displaystyle
(p\cutc)^{-\halfff 1}\dff(S\smallsetminus c)\ttoo S\smallsetminus c$\oss 
and a double covering map\oss 
\[
\quad
({p\cutc})^{-\halfff 1}\dff(c)\off\toto\off c\fff.
\]
We will treat the diffeomorphism\oss
$\displaystyle
(p\cutc)^{-1}\dff(S\smallsetminus c)\ttoo S\smallsetminus c$\oss 
as an identification.\oss
If\qss ${p\cutc}$\qss is injective on a component\dss $Q$\dss of\oss $S\cutc$\snsp,\qff\oss 
then we will treat  as an identification also the induced map\dss 
$\displaystyle
Q\ttoo p\cutc\dff(Q)$\dnsp,\oss 
and will treat\dss $Q$\dss as a subsurface of\dss $S$\dnsp.\oss 

Since\dss $S$\dss is orientable,\oss the covering\oss
$\displaystyle
({p\cutc})^{-\halfff 1}\dff(c)\ttoo c$\oss is actually trivial\halfff.\oss
Namely,\oss for every component\dss $C$\dss of\dss $c$\dss its preimage\qss 
$(p\cutc)^{-1}\dff(C)$\qss 
consists of two boundary circles of\qss $S\cutc$\dnsp,\oss
and each of these circles is mapped by\dss $p\cutc$\dss diffeomorphically onto\dss $C$\dfdot 

Any diffeomorphism\dss $F\fff\colon S\ttoo S$\dss such that\qss 
$F\fff(c)\qff =\qff c$\qss induces a diffeomorphism 
\[
\quad
F\cutc\dff\colon\dff S\cutc\ttoo S\cutc\dff.
\] 
If\pss $F\cutc$\qss leaves a component\dss $Q$\dss 
of\qss $S\cutc$\qss 
invariant\halfff,\oss 
then\qss $F\cutc$\qss induces a diffeomorphism\vspace*{-0.5\medskipamount}
\[
\quad
F_Q\dff\colon Q\ttoo Q\fff,
\]

\vspace*{-1.25\bigskipamount}
called the\pss \emph{restriction}\pss of\qss $F$\qss to\qss $Q$\dnsp.\oss

\myitpar{Systems of circles.}\qss A one-dimensional closed submanifold\dss $c$\dss 
of a surface\dss $S$\dss is called a\qss 
\emph{system of circles on}\dss $S$\qss if the components of\dss $c$\dss are all
non-trivial circles on\dss $S$\dss and are pair-wise non-isotopic.\oss
As usual,\oss
we will denote by\dss $\pi_{\dff 0\fff}(c)$\dss the set of components of\dss $c$\dfdot
If\dss $c$\dss is a system of circles,\qss 
then the elements of\dss $\pi_{\dff 0\fff}(c)$\dss 
are pair-wise disjoint and pair-wise non-isotopic circles,\oss
and if\dss $U$\dss is a set of pair-wise disjoint and non-isotopic circles,\qss
then the union of these circles is a system of circles.\oss

\myitpar{Reduction systems.} A system of circles\dss $c$\dss on\dss $S$\dss 
is called a\qss \emph{reduction system for a diffeomorphism}\qss $F\fff\colon S\ttoo S$\qss
if\qss $F\fff(c)\qff =\qff c$\snsp.\oss
A system of circles\dss $c$\dss on\dss $S$\dss 
is called a\qss \emph{reduction system for an element}\qss $f\dff\in\dff \mms$\qss
if\dss $c$\dss is a reduction system for some diffeomorphism in the isotopy class\qss $f$\snsp,\oss
i.e.\qss if the isotopy class\qss $f$\qss can be represented by a diffeomorphism\qss
$F\fff\colon S\ttoo S$\qss such that\qss $F\fff(c)\qff =\qff c$\dnsp.\oss

\myitpar{Reducible and pseudo-Anosov elements.}
A non-trivial\qss element\qss $f\dff\in \mmod(S)$\qss is said to be\qss \emph{reducible}\qss 
if\dss there exists a non-empty reduction system for\qss $f$\halfff\dnsp,\oss
and\qss \emph{irreducible}\qss otherwise.\qss
An irreducible element of infinite order is called a\qss \emph{pseudo-Anosov element}.\oss
This is an easy\nsp,\oss but hardly enlightening{\halfff} way to define pseudo-Anosov elements.\oss
The original definition of\qss Thurston\qss \cite{thurston}\qss
looks more like a theory than a short definition.\oss

\mysection{Dehn\qss twists}{dehn-twist}

\myitpar{Twist diffeomorphisms of an annulus.} Let $A$ be an annulus,\oss
i.e. a surface diffeomorphic to $S^1\dnsp\times\nsp [\halfff 0,\dff 1]$\dfdot
As is well known,\oss the group of diffeomorphisms of\dss $A$\halfff\dfcom
fixed in a neighborhood of\dss $\partial A$\dss and considered up to isotopies 
fixed in a neighborhood\dss $\partial A$\halfff\dfcom
is an infinite cyclic group.\oss
A diffeomorphism of\dss $A$\dss is called a\qss \emph{twist diffeomorphism}\qss of\dss $A$\dss
if it is fixed in a neighborhood of\dss $\partial A$\dss and its isotopy class 
is a generator of this group.\oss 

Fixing an orientation of\dss $A$\dss allows to choose a 
preferred generator of this group.\oss
A twist diffeomorphism of\dss $A\ttoo A$\dss 
is called a\qss \emph{left twist diffeomorphism}\qss if its isotopy class is 
the preferred generator\halfff,\oss and a\qss{right twist diffeomorphism}\qss otherwise.\oss
Up to an isotopy,\oss the right twist diffeomorphism is the inverse of the left one.

\myitpar{Twist diffeomorphisms of a surface.} 
If\dss $A$\dss is annulus contained in\dss $S$\dss as a subsurface,\oss
then any twist diffeomorphism of $A$ can be extended by 
the identity to a diffeomorphism of\qss $S$\dnsp.\oss
Such extensions are called\qss \emph{twist diffeomorphisms}\qss of\qss $S$\dnsp.\oss

If\qss $S$\qss is oriented,\oss then the extension of a left\qss
(respectively,\oss right)\qss twist diffeomorphism of an annulus in\qss $S$\dss 
with the induced orientation is called a\pss
\emph{left}\pss (respectively,\oss \emph{right})\qss 
twist diffeomorphism of\qss $S$\dnsp.\oss

\myitpar{Dehn twists.}\qss Suppose that\dss $S$\dss is oriented.\oss
Let\qss $A$\qss be an annulus contained in\qss $S$\dnsp,\oss
and let\qss $C$\qss be a circle in\qss $S$\qss contained in\qss $A$\qss
as a deformation retract\halfff.\oss
Such a circle\qss $C$\qss is unique up to isotopy.\oss

It turns out that isotopy class of a left\qss (respectively,\oss right)\qss 
twist diffeomorphism of\qss $S$\qss
obtained by extension of a left\qss (respectively,\oss right)\qss 
twist diffeomorphism of\dss $A$\dss 
depends only on the isotopy class of\qss $C$\qss in\qss $S$\dnsp.\oss
In particular\halfff,\oss it is uniquely determined by\qss $C$\dnsp.\oss
This isotopy class is called the\qss \emph{left}\pss
(respectively,\oss \emph{right})\pss 
\emph{Dehn twist}\qss of\qss $S$\qss about the circle\qss $C$\dnsp.\qff\oss
It is an element of the group\qss $\mms$\dnsp.\oss

The left Dehn twist about a circle\qss $C$\qss is denoted by\pss $t_{\dff C}$\nsp.\oss
The right Dehn twist about\qss $C$\qss is the inverse of the left one
and hence is equal to\pss $t_{\dff C}^{\dff -\halfff 1}$\nsp.\oss

The Dehn twist\qss $t_{\dff C}$\qss
is equal to\qss $1\dff\in\dff \mms$\qss if and only if the circle\qss $C$\qss is trivial.\oss
Let\qss $G$\qss be a diffeomorphism of\qss $S$\qss and 
let\qss $g\dff\in\dff \mms$\qss be its isotopy class.\oss
Then\oss
\[
\quad
g\qff t_{\dff C}\qff g^{\fff -\dff 1}\qff =\off t_{\dff G\fff(C)}\fff.
\]
In particular\halfff,\oss if\qss $G\fff(C)\qff =\qff C$\snsp,\oss
or\qss if\qss $G\fff(C)$\qss is isotopic to\qss $C$\snsp,\oss
then\qss $t_{\dff C}$\qss and\qss $g$\qss commute.\oss

\mysection{Action\qss of\qss Dehn\qss twists\qss on\qss homology}{dehn-torelli}

\myitpar{The explicit formula.}\qss 
\emph{From now on we will assume that the surface\qss $S$\qss is closed and oriented.\oss}
The orientation of  $S$ allows to define a
skew-symmetric pairing on $H_1(S)$\dfcom known as the\qss \emph{intersection pairing}.\oss
It is denoted by
\[
\quad
(a\fff,\qff b)\off \longmapsto\off \langle a\fff,\qff b\rangle\dff.
\]

\vspace*{-0.25\medskipamount}
\vspace*{-\bigskipamount}
Let $C$ be a circle on\dss $S$\dfdot
Let us orient $C$ and denote by $\hclass{C}$
the image of the fundamental class of\dss $C$\dss in\qss $H_{\fff 1}(S)$\dfdot
Then\qss $t_{\dff C}$\qss acts on\qss $H_{\fff 1}(S)$\qss by the formula\oss
\begin{equation}
\label{twist-homology}
\quad
\bigl(t_{\dff C}\bigr)_*\fff(a)
\off =\off 
a\qff +\qff  \langle\fff a\fff,\qff \hclass{C}\dff\rangle\dff\hclass{C}\dff.
\end{equation}
Changing the orientation of $C$ replaces $\hclass{C}$ by $- \fff\hclass{C}$
and hence does not change the right hand side of\qss (\ref{twist-homology}).\oss
The powers of\qss $t_{\fff C}$\qss act on\qss $H_{\fff 1}(S)$\qss by the formula\oss 
\begin{equation}
\label{right-twist-homology}
\quad
\bigl(t_{\dff C}^{\dff m_{\phantom{O}}}\hspace*{-0.5em}\bigr)_*\fff(a)\off 
=\off a\qff +\qff m\dff \langle\fff a\fff,\qff \hclass{C}\dff\rangle\dff\hclass{C}\dff.
\end{equation}
If\qss $\hclass{C}\qff \neq\qff 0$\dnsp,\qff\oss 
then\qss $\langle a\fff,\qff \hclass{C}\dff\rangle\qff \neq\qff 0$\qss
for some\qss $a\dff\in\dff H_1(S)$\dnsp,\oss
and hence\oss 
\[
\quad
a\qff +\qff  \langle\fff a\fff,\qff \hclass{C}\dff\rangle\dff[C]\off \neq\off a\dff.
\]
On the other hand,\oss if\qss $\hclass{C}\qff =\qff 0$\dnsp,\qff\oss 
then\oss
$\displaystyle 
a\qff +\qff  \langle\fff a\fff,\qff \hclass{C}\dff\rangle\dff[C]
\off =\off a$\dnsp.\qff\oss
Therefore,\oss $t_{\fff C}\dff\in\dff\tors$\qss
if and only if\oss 
$\displaystyle
\hclass{C}\off =\off 0$\snsp.\oss
i.e.\qss if\dss and\sss only\sss if\qss $C$\qss is a separating circle.\oss

\myitpar{Dehn--Johnson twists and multi-twists.} Let\qss $C\fff,\pff D$\qss 
be a pair of disjoint circles on\qss $S$\qss 
such that the union\qss $C\dff\cup\dff D$\qss is equal to the boundary\qss
$\partial Q$\qss of some subsurface\qss $Q$\qss of\qss $S$\dnsp.\oss
Such pairs of circles are called\pss \emph{bounding pairs}.\qff\oss
The circles\dss $C$\dss and\dss $D$\dss can be oriented 
in such a way that\qss $\hclass{C}\off =\off \hclass{D}$\dnsp.\qff\oss
Therefore\qss (\ref{twist-homology})\qss 
implies that the maps
\[
\quad
\bigl(t_{\fff C}\bigr)_*\qff,\off\qff \bigl(t_{\fff D}\bigr)_*\qff\colon\qff H_1(S)\qff\ttoo\qff H_1(S)
\]
are equal,\oss and hence\oss
$\displaystyle
t_{\fff C}\trf t_{\fff D}^{\fff -\halfff 1}\fff,\off\off 
t_{\fff D}\trf t_{\fff C}^{\fff -\halfff 1}\qff \in\off \tors$\dnsp.\qff\oss
Both these elements of\qss $\tors$\qss are called the\pss 
\emph{Dehn--Johnson twists}\pss about the bounding pair\qss $C\fff,\pff D$\dnsp.\oss 

Let\dss $c$\dss be a one-dimensional closed submanifold of\qss $S$\dfdot
A\qss \emph{Dehn multi-twist about}\qss about\dss $c$\dss 
is defined as a product\qss $t$\qss of the form
\begin{equation}
\label{multi-twist-definition}
\quad
t\off =\off \prod\nolimits_{\fff O}\qff t_{\fff O}^{\dff m_{\dff O}},
\end{equation}
where\qss $O$\qss runs over all components of\qss $c$\snsp,\oss
and\qss $m_{\dff O}$\qss are integers.\qss

\myitpar{A siren song.} Suppose that\qss $Q$\qss 
is a subsurface of\qss $S$\qss and\qss $c\qff =\qff \partial Q$\dnsp.\qff\oss
The orientation of\qss $S$\qss defines an orientation of\qss $Q$\dnsp,\oss
which,\oss in turn,\oss defines an orientation of\pss $c\off =\off \partial Q$\dnsp.
Let\qss $t$\qss be defined by\oss (\ref{multi-twist-definition}).\oss
If\qss $c$\qss has\qss $\leq 2$\qss components,\oss
then\oss $t\dff\in\dff\tors$\oss if and only if
\begin{equation}
\label{homology-class}
\quad
\sum\nolimits_{\fff O}\qff m_{\dff O}\dff \hclass{O}\off =\off 0\dff,
\end{equation}
where\qss $O$\qss runs over components of\dss $c$\dnsp,\oss
and these components are considered with orientations 
induced from\qss $c\off =\off \partial Q$\dnsp.\qff\oss
It is tempting to believe that this is true in general.\oss

In fact\halfff,\oss this is very far from being true.\oss
The following example,\oss
explained to the author in few words by\qss R.\qss Hain\dss in\dss 1992,\oss 
illustrates the reason behind this,\oss
and Theorem\qss \ref{multi-twists-in-t}\qss below provides
a simple necessary and sufficient condition for\qss $t\dff\in\dff \tors$\dnsp.\oss

It is much more difficult to find reasons for such a belief.\oss

\myitpar{Example.}\qss
Suppose that both\qss $Q$\qss and the complementary subsurface\qss $\ccomp Q$\qss 
are connected and that\qss $c\qff =\qff \partial Q$\qss 
consists of three circles\qss $C\fff,\pff D\fff,\pff E$\snsp.\oss
Let us orient these circles as boundary components of\qss $Q$\dnsp.\qff\oss
Then the only relation between\oss $\hclass{C}$\snsp,\oss 
$\hclass{D}$\snsp,\oss $\hclass{E}$\oss is
\[
\quad
\hclass{C}\qff +\qff \hclass{D}\qff +\qff \hclass{E}\off =\off 0\dff.
\]
In particular\halfff,\qss the classes\oss $\hclass{C}$\snsp,\oss 
$\hclass{D}$\snsp,\oss $\hclass{E}$\oss
are non-zero and any two of them are linearly independent\halfff.\oss
Let\dss $A$\dss be a circle in\dss $S$\dss disjoint from\dss $E$\dss and 
intersecting each of the circles\qss $C\fff,\pff D$\qss transversely at one point\halfff.\oss
Let us orient\dss $A$\dss and let\qss $a\qff =\qff \hclass{A}$\dfdot
By the choice of\dss $A$\dss the intersection number\qss 
$\langle\fff a\fff,\qff \hclass{E}\dff\rangle\qff =\qff 0$\snsp,\oss
one of the intersection numbers\oss 
$\langle\fff a\fff,\qff \hclass{C}\dff\rangle$\dnsp,\oss
$\langle\fff a\fff,\qff \hclass{D}\dff\rangle$\oss is equal to\dss $1$\snsp,\pss 
and the other is equal to\qss $-\halfff 1$\dfdot
We may assume that\pss 
$\langle\fff a\fff,\qff \hclass{C}\dff\rangle\qff =\qff 1$\pss 
and\pss
$\langle\fff a\fff,\qff \hclass{D}\dff\rangle\qff =\qff -\halfff 1$\dnsp.\qff\oss
Let\pss
$m_{\dff C}$\snsp,\oss $m_{\dff D}$\snsp,\oss $m_{\dff E}\pff\in\pff \zzz$\oss 
and\oss
\[
\quad
t\off =\off t_C^{\dff m_{\dff C}}\pff t_D^{\dff m_{\dff D}}\pff t_C^{\dff m_{\dff D}}\fff.
\]
Then\hspace*{0.7em}
$\displaystyle
t_{\fff *}\dff(a)\off =\off a\qff +\qff m_{\dff C}\dff \langle\fff a\fff,\qff 
\hclass{C}\dff\rangle\dff \hclass{C}\qff
+\qff  m_{\dff D}\dff \langle\fff a\fff,\qff \hclass{D}\dff\rangle\dff \hclass{D}\qff
+\qff  m_{\dff E}\dff \langle\fff a\fff,\qff \hclass{E}\dff\rangle\dff \hclass{E}\off
$
\[
\quad
\phantom{t_{\fff *}\dff(a)\off }
=\off a\qff +\qff m_{\dff C}\dff \hclass{C}\qff -\qff m_{\dff D}\dff \hclass{D}.
\]

\vspace*{-\medskipamount}
It follows that if\pss $t\dff\in\dff\tors$\dfcom
then\qss $m_{\dff C}\qff =\qff m_{\dff D}\qff =\qff 0$\dfdot
A similar argument using\qss $D\fff,\pff E$\qss instead of\qss $C\fff,\pff D$\qss shows that
if\pss $t\dff\in\dff\tors$\dfcom then also\qss $m_{\dff D}\qff =\qff m_{\dff E}\qff =\qff 0$\dfdot
It follows that\pss $t\dff\in\dff\tors$\pss 
if and only if\oss $m_{\dff C}\qff =\qff m_{\dff D}\qff =\qff m_{\dff E}\qff =\qff 0$\dnsp,\off\oss
i.e.\qss if and only if\oss $t\qff =\qff 1$\dfdot

\mysection{Dehn\qss multi-twists\qss in\qss Torelli\qss groups}{multi-twists-torelli}

\myitpar{Homology equivalence.} 
Two circles\qss $C\fff,\pff C'$\qss in\qss $S$\qss are said to be\qss 
\emph{homology equivalent}\qss if they can be oriented in such way that their
homology classes in\qss $H_{\fff 1}(S)$\qss are equal,\oss 
\[
\hclass{C}\qff =\qff \hclass{C'}\dff.
\]
Changing the orientation of a circle\qss $C$\qss replaces\qss 
$\hclass{C}$\qss by\qss $-\qff \hclass{C}$\dnsp.\oss 
It follows that the homology equivalence is indeed an equivalence relation.\oss
If $C$ is a separating circle, then $\hclass{C}\qff =\qff 0$\dnsp,\oss 
and hence all separating circles are homology equivalent.

\myitpar{Necklaces.} Let\qss $s$\qss be a one-dimensional closed submanifold of\qss $S$\dnsp.\oss
The homology equivalence induces an equivalence relation 
on the set\qss $\pi_{\dff 0\fff}(s)$\qss of components of\dss $s$\dnsp.\oss
A\qss \emph{necklace}\qss of\qss $s$\qss is an equivalence class of this equivalence
relation containing a non-separating circle.\oss
If a necklace consists of more than $1$ circle,\oss
then the union of all circles in it is a\qss \emph{BP-necklace}\qss
in the sense of Appendix.\oss
See Corollary\qss \ref{BP-and-necklaces}.\oss
The description of BP-necklaces in Appendix is
the motivation behind the term\qss \emph{necklace}.\oss

\myitpar{Graph associated with a submanifold.} Let\qss $c$\qss be a 
one-dimensional closed submanifold of\qss $S$\dnsp.\oss
One can associate with the pair\qss $(S\fff,\pff c)$\qss
a graph\qss $\mathfrak{G}(S\fff,\pff c)$\qss as follows.\oss
Its set of vertices is the set\qss $\pi_{\dff 0\fff}(S\cutc)$\qss 
of components of the cut surface\qss $S\cutc$\dnsp,\oss
and its set of edges is the set\qss $\pi_{\dff 0\fff}(c)$\qss of components of\dss $c$\dnsp.\oss
Every component\qss $C$\qss of\qss $c$\qss is the image under the canonical map\qss
$p\qff =\qff p\cutc\dff\colon S\cutc\ttoo S$\qss of two components of\qss $\partial S\cutc$\dnsp.\oss
Let\qss $C_{\fff 1}\fff,\pff C_{\fff 2}$\qss be these components,\oss
and let\qss $Q_1\fff,\pff Q_2$\qss be the components of\qss $S\cutc$\qss containing\qss 
$C_{\fff 1}\fff,\pff C_{\fff 2}$\qss respectively\qss
(it may happen that\qss $Q_1\qff =\qff Q_2$\nsp).\oss
The component\qss $C$\qss considered as an edge of\qss 
$\mathfrak{G}(S\fff,\pff c)$\qss connects 
the components\qss $Q_1\fff,\pff Q_2$\qss considered as vertices.\oss
If\qss $Q_1\qff =\qff Q_2$\dnsp,\oss then the edge\dss $C$\dss is a loop,\oss
and if the intersection of the images\qss $p(Q_1)\fff,\pff p(Q_2)$\qss
consists of several components of\dss $c$\dnsp,\oss
then the vertices\qss $Q_1\fff,\pff Q_2$\qss are connected by several edges.\oss

\myitpar{Separating edges.} Let\qss $a\fff,\pff b$\qss be two vertices of 
a connected graph\qss $\mathfrak{G}$\dnsp.\oss
A set\qss $T$\qss of edges of\qss $\mathfrak{G}$\qss is\qss
\emph{separating edge set}\qss for\qss $a$\qss and\qss $b$\qss if every path in\qss
$\mathfrak{G}$\qss contains at least one edge from\qss $T$\dnsp.\oss
The vertices\qss $a\fff,\pff b$\qss are said to be\qss $\tau$\dnsp-edge separated
if\qss $\tau$\qss is the minimal number of elements in 
a separating edge set for\qss $a$\qss and\qss $b$\dnsp.\oss

\mypar{Theorem.}{edge-separation} \emph{If two vertices\qss $a\fff,\pff b$\qss
of a connected graph\qss $\mathfrak{G}$\qss are\qss $\tau$\dnsp-edge separated,\oss
then there exists\qss $\tau$\qss simple paths connecting\qss $a$\qss with\qss $b$\qss
and having pair-wise disjoint sets of edges.\oss}

\proof\qss See,\oss for example,\oss Theorem 12.3.1\qss in the classical book of\qss 
O.\qss Ore\qss \cite{o}.\oss  \eproof

\mypar{Corollary.}{one-edge-separation} \emph{If two vertices\qss $a\fff,\pff b$\qss
of a connected graph\qss $\mathfrak{G}$\qss remain connected by a path after removing
any edge of\qss $\mathfrak{G}$\dnsp,\oss
then there exists at least two\qss $2$\qss simple paths connecting\qss $a$\qss with\qss $b$\qss
and having disjoint sets of edges.\oss}

\myitpar{Remark.} We will use only Corollary\qss \ref{one-edge-separation}.\oss
But Theorem\qss \ref{edge-separation}\qss is so beautiful that the author could
not resist including it\halfff.\oss
It is an edge-separation version of vertex-separation results
of\dss K. Menger\qss \cite{menger}\qss and H. Whitney\qss \cite{w}.\oss 
It is worth to note that\dss K.\qss Menger's paper\qss \cite{menger}\qss is a paper in topology,\oss
and\dss H.\qss Whitney is better known not for his seminal contributions to graph theory,\oss 
but as one of the creators of differential topology.\oss

The idea to apply Corollary\qss \ref{one-edge-separation}\qss 
to Dehn multi-twists in Torelli groups is due to W. Vautaw\qss \cite{v}.\oss
Unfortunately,\oss the application of this result is hidden deep inside of a
technical argument in\qss \cite{v}\halfff,\oss 
and no references related to this result are provided.\oss

For a textbook exposition the vertex-separation version of\dss 
Corollary\qss \ref{one-edge-separation},\oss known as Whitney's theorem,\oss
see\qss \cite{bm},\oss Theorem\qss 3.2,\oss or\qss
\cite{br},\oss Theorem\qss 3.3.7.\oss 
The Corollary\qss \ref{one-edge-separation}\qss 
itself is the Exercise\qss 3.9\qss in\qss \cite{bm}.\oss

\mypar{Lemma.}{two-circles} \emph{Suppose that\dss $c$\dss is a one-dimensional
submanifold of\qss $S$\dss such that all components of\dss $c$\dss are non-separating
and no two components are homology equivalent\halfff.\oss
Let\dss $D$\dss be a component of\dss $c$\snsp.\qff\oss
Then there exists two circles\qss $A\fff,\pff B$\qss on\qss $S$\qss
intersecting each component of\qss $c$\qss transversely in no more than one point
and such that\dss $D$\dss intersects both circles\qss $A\fff,\pff B$\qss and no other 
component of\dss $c$\dss does.}

\proof\qss Let us consider the graph\qss $\mathfrak{G}(S\fff,\pff c)$\dnsp.\oss
Let\qss $Q_1\fff,\pff Q_2$\qss be the components of\qss $S\cutc$\qss connected by
the edge\qss $D$\qss of\qss $\mathfrak{G}(S\fff,\pff c)$\dnsp.\oss
If\qss $Q_1\qff =\qff Q_2$\nsp,\oss then there is a circle in\qss $S$\qss
intersecting\qss $D$\qss transversely in one point and disjoint from all
other components of\dss $c$\dnsp.\oss
In this case we can take this circle as both\qss $A$\qss and\qss $B$\dnsp.\oss

Suppose now that\qss $Q_1\qff \neq\qff Q_2$\nsp.\oss
Let\qss $\mathfrak{G}$\qss be the result of removing 
the edge\dss $D$\dss from the graph\qss $\mathfrak{G}(S\fff,\pff c)$\dnsp.\oss
Since the circle\dss $D$\dss is non-separating,\oss
the graph\qss $\mathfrak{G}$\qss is connected.\oss
If after removing an edge\dss $C$\dss of\qss $\mathfrak{G}$\qss
the vertices\qss $Q_1\fff,\pff Q_2$\qss are not connected,\oss
then\qss $C$\qss and\qss $D$\qss together separate\qss $S$\dnsp,\oss
and hence\dss $C$\dss and\dss $D$\dss together bound a subsurface of\qss $S$\dnsp.\oss
In this case the circles\dss $C$\dss and\dss $D$\dss are homology equivalent\halfff,\oss
contrary to the assumption.\oss
Therefore the vertices\qss $Q_1\fff,\pff Q_2$\qss of\qss $\mathfrak{G}$\qss remain
connected after removing any edge of\qss $\mathfrak{G}$\dnsp.\oss

Let us choose some points\qss $x_1\fff,\pff x_2$\qss in the interior of surfaces\qss
$Q_1\fff,\pff Q_2$\qss respectively.\oss
Let\qss $J$\qss be an arc in\qss $S$\qss connecting\qss $x_1$\qss with\qss $x_2$\nsp,\oss
intersecting\qss $D$\qss transversely at one point\halfff,\oss
and disjoint from all other components of\dss $c$\dnsp.\oss

A path in\qss $\mathfrak{G}$\qss connecting\qss $Q_1$\qss with\qss $Q_2$\qss
can be\qss ``realized''\qss by an arc in\qss $S$\qss connecting\qss $x_1$\qss with\qss $x_2$\nsp,\oss
contained in the union
\[
\quad
\bigcup\nolimits_{\fff Q}\dff p(Q)\fff,
\]
where\qss $Q$\qss runs over the set of vertices of this path\halfff,\oss
intersecting only those components of\dss $c$\dss which are the edges of this path\halfff,\oss
and intersecting each of these components transversely at one point\halfff.\oss 

By Corollary\qss \ref{one-edge-separation},\oss there exist two paths in\qss $\mathfrak{G}$\qss
connecting\qss $Q_1$\qss with\qss $Q_2$\qss and having disjoint sets of edges.\oss
Therefore,\oss there exist two arcs in\qss $S$\qss connecting\qss $x_1$\qss with\qss $x_2$\nsp,\oss
disjoint from\qss $D$\dnsp,\oss and such that no other component of\dss $c$\dss intersects both of these arcs.\oss
By taking the unions of these two arcs with the arc\dss $J$\dss we get circles\qss $A\fff,\pff B$\qss
with the required properties.\oss  \eproof

\myitpar{Difference maps.} The\qss \emph{difference map}\qss 
of a diffeomorphism\qss $G\dff\colon S\ttoo S$\qss is the map
\[
\quad
\Delta_{\dff G}\dff\colon H_{\fff 1}(S)\ttoo H_{\fff 1}(S)
\]
defined by\oss 
$\displaystyle
\Delta_{\dff G}\off =\off G_{\fff *}(a)\qff -\qff a$\snsp.\oss
Clearly,\oss $\Delta_{\dff G}$\qss depends only on the isotopy class\qss $g$\qss of\qss $G$\qss
and hence may be denoted by\dss $\Delta_{\dff g}$\dnsp.\oss
The isotopy class\qss $g$\qss belongs to\qss $\tors$\qss if and only if\qss 
the difference map\qss $\Delta_{\dff g}$\qss is equal to\qss $0$\snsp.\oss

Let\qss $c$\qss be a one-dimensional closed submanifold of\qss $S$\dnsp,\oss
and let\oss
\[
\quad
u\off =\off
\prod\nolimits_{\fff C}\qff t_{\fff C}^{\dff n_{\fff C}}\qff\in\qff \tors\dff,
\]
where\qss $C$\qss runs over all components of\dss $c$\dss and\qss
$n_{\fff C}$\qss are integers,\oss be a Dehn multi-twist about\dss $c$\dnsp.\oss
The formula\qss (\ref{twist-homology})\qss allows to compute
the difference map\qss $\Delta_{\dff u}$\nsp.\oss Namely,\oss
\begin{equation}
\label{difference-multi-twist}
\quad
\Delta_{\dff u}\dff(a)
\off =\off
\sum\nolimits_{\fff C}\qff n_{\fff C}\trf \langle\fff a\fff,\qff \hclass{C}\dff\rangle\dff\hclass{C}\dff,
\end{equation}
where\qss $C$\qss runs over all components of\dss $c$\dnsp.\oss

\mypar{Theorem.}{multi-twists-in-t} \emph{Let\qss $s$\qss 
be a one-dimensional closed submanifold of\qss $S$\dnsp,\oss
and let\oss $t$\pss be a Dehn multi-twist 
about\qss $s$\dnsp,\oss i.e.\oss 
\begin{equation*}
\quad
t\off =\off \prod\nolimits_{\fff O}\qff t_{\fff O}^{\dff m_{\fff O}},
\end{equation*}
where\qss $O$\qss runs over all components of\qss $s$\dnsp,\oss
and\qss $m_{\fff O}\dff\in\dff \zzz$\nsp.\qff\oss
Then\qss $t\qff\in\qff \tors$\qss if and only if}
\[
\quad
\sum\nolimits_{\fff O\dff\in\dff N}\qff m_{\fff O}\off =\off 0
\]

\vspace*{-0.7\bigskipamount}
\emph{for all necklaces\qss $N$\qss of\pss $s$\dnsp.\oss}

\proof\qss Since all Dehn twists about separating circles belong to\qss $\tors$\dnsp,\oss
we may assume that\qss $s$\qss has no separating components.\oss
Let us select a circle from every necklace of\dss $s$\dnsp.\oss
Let\qss $c$\trs be the union of all selected circles.\oss
For every selected circle\qss $C$\qss let\oss
\[
\quad
n_{\fff C}\off =\off \sum\nolimits_{O}\qff m_{\fff O}\dff,
\]
where the sum is taken over all circles in the necklace containing\qss $C$\dnsp.\oss
Since Dehn twists about homology equivalent circles 
induce the same automorphism of\qss $H_{\fff 1}(S)$\dnsp,\oss
the Dehn multi-twist\qss $t$\qss belongs to\qss $\tors$\qss if an only if
\[
\quad
u\off =\off
\prod\nolimits_{\fff C}\qff t_{\fff C}^{\dff n_{\fff C}}\qff\in\qff \tors\dff,
\]
where\qss $C$\qss runs over all components of\dss $c$\dnsp.\oss 

We need to prove that\qss $u\dff\in\dff \tors$\qss if
and only if\qss $n_{\fff C}\qff =\qff 0$\qss for all components\qss $C$\qss of\qss $c$\dnsp.\qff\oss
If\dss all\qss $n_{\fff C}\qff =\qff 0$\dnsp,\oss 
then\qss $u\qff =\qff 1\dff\in\dff \tors$\dnsp.\oss
This proves the\qss \emph{``if''}\qss part\halfff.\oss 

Suppose now that\qss $u\dff\in\dff \tors$\dnsp.\oss
By the choice of\dss $c$\dnsp,\oss all components of\dss $c$\dss are non-separating
and no two of them are homology equivalent\halfff.\oss
Let\qss $D$\qss be a component of\dss $c$\dnsp,\oss
and let\qss $A\fff,\pff B$\qss be two circles provided by Lemma\qss \ref{two-circles}.\qff\oss
Let us orient circles\qss $A\fff,\pff B\fff,\pff D$\qss and consider their homology
classes\qff\oss 
$\displaystyle
a\qff =\qff \hclass{A}\fff,\quad 
b\qff =\qff \hclass{B}\fff,\quad 
\hclass{D}$\nsp.\oss
We may assume that\oss
$\displaystyle
\langle\fff a\fff,\qff \hclass{D}\trf\rangle
\off =\off 
\langle\fff a\fff,\qff \hclass{D}\trf\rangle
\off =\off 
1$\nsp.\oss
The formila\qss (\ref{difference-multi-twist})\qss implies that
\[
\quad
\langle\dff b\fff,\qff \Delta_{\dff u}\dff(a)\dff\rangle
\off =\off
\sum\nolimits_{\fff C}\qff n_{\fff C}\trf \langle\fff a\fff,\qff \hclass{C}\dff\rangle\trf
\langle\dff b\fff,\qff \hclass{C}\dff\rangle\dff\dff,
\]
where\qss $C$\qss runs over all components of\dss $c$\dnsp.\oss
By the choice of circles\qss $A\fff,\pff B$\dnsp,\qff\oss for all components\qss $C\qff \neq\qff D$\qss
either\oss $\langle\fff a\fff,\qff \hclass{C}\dff\rangle\off =\off 0$\dnsp,\qff\oss or\oss
$\langle\fff b\fff,\qff \hclass{C}\dff\rangle\off =\off 0$\nsp.\oss
It follows that\oss
\[
\quad
\langle\dff b\fff,\qff \Delta_{\dff u}\dff(a)\dff\rangle
\off =\off\dff n_{\fff D}\fff.
\]
On the other hand,\oss $\Delta_{\dff u}\dff(a)\qff =\qff 0$\qss because\qss $u\dff\in\dff\tors$\dnsp.\oss
It follows that\qss $n_{\fff D}\qff =\qff 0$\dnsp.\oss
Since the component\qss $D$\qss of\dss $c$\dss was arbitrary,\oss
this proves the\qss \emph{``only if''}\qss part of the theorem.\oss  \eproof

\mysection{The\qss rank\qss of\qss multi-twist\qss 
subgroups\qss of\qss Torelli\qss groups}{multi-twist-subgroups}

\vspace*{\bigskipamount}
The rest of the paper is focused on the estimates of the rank of various
abelian subgroups of\qss $\tors$\dnsp.\oss
We denote the rank of an abelian group\qss $\mathfrak{A}$\qss
by\qss $\rank \mathfrak{A}$\oss its rank\halfff.\oss

\myitpar{Dehn multi-twist subgroups.} Let\dss $s$\dss 
be a one-dimensional closed submanifold of\qss $S$\dnsp.\qff\oss
By\qss $\mathcal{T}\fff(s)$\qss we will denote the subgroup of\qss $\mms$\qss
generated by Dehn twists about components of\dss $s$\dnsp,\oss
i.e.\qss the group of Dehn multi-twists about\dss $s$\dnsp.\oss
If\dss $s$\dss is a system of circles,\oss then\qss $\mathcal{T}\fff(s)$\qss
is well known to be a free abelian group having\qss (say,\oss left\halfff)\qss Dehn twists about
components of\dss $s$\dss as free generators.\oss

\mypar{Theorem.}{multi-twists-subgroups} \emph{Let\dss $s$\dss 
be a system of circles on\qss $S$\snsp.\qff\oss
Then\oss
\[
\quad
\rankf\dff \mathcal{T}\fff(s)\dff\cap\dff \tors \off =\off N\qff -\qff n\fff,
\]
where\qss $N$\qss is the number of components of\pss $s$\snsp,\oss
and\qss $n$\qss is the number of necklaces of\pss $s$\snsp.\oss}

\proof\qss Since\oss $\mathcal{T}\fff(s)$\qss is a 
free abelian group of rank\qss $N$\qss
freely generated by Dehn twists about components of\qss $s$\snsp,\oss
this theorem immediately follows from Theorem\qss \ref{multi-twists-in-t}.\oss  \eproof

\begin{figure}[t]
\includegraphics[width=0.96\textwidth]{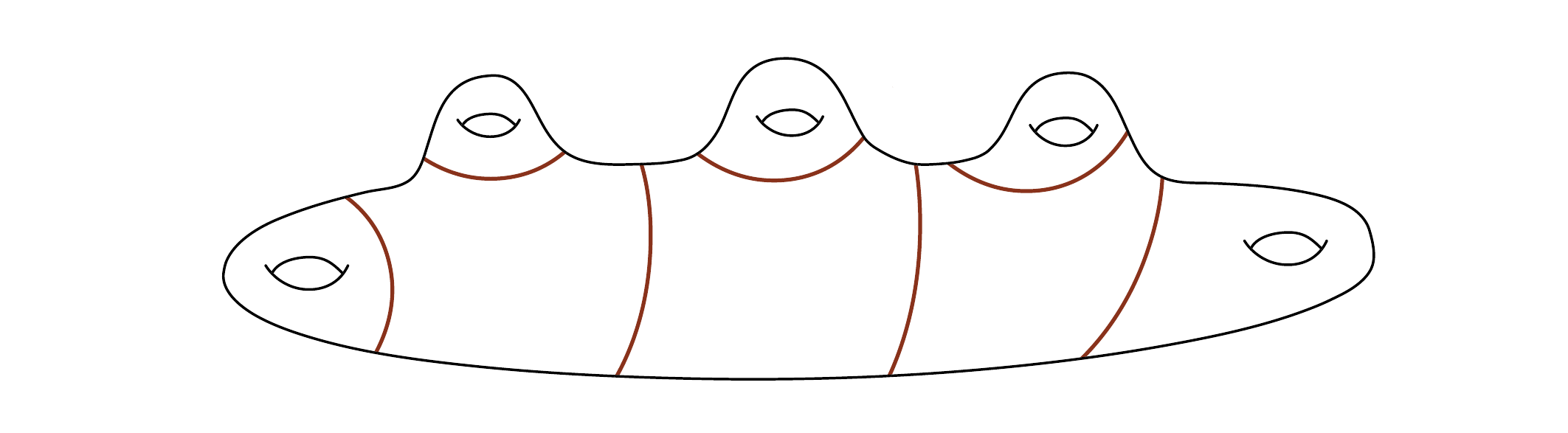}
\caption{A system of\dss $2\dff\gen\qff -\qff 3$\dss separating circles.}
\label{separating-circles}
\end{figure}

\begin{figure}[h]
\includegraphics[width=0.96\textwidth]{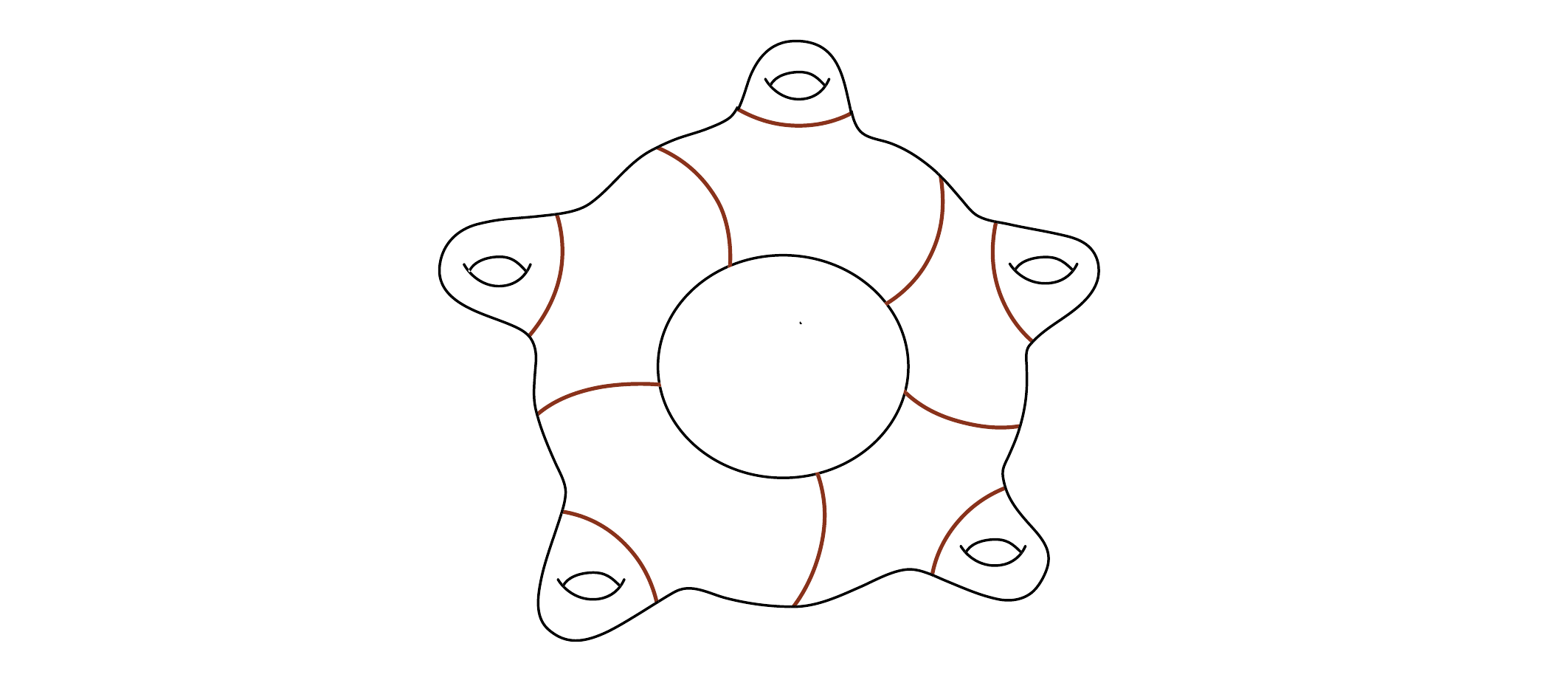}
\caption{A system of\dss $\gen\qff -\qff 1$\dss separating and\dss $\gen\qff -\qff 1$\dss 
non-separating circles.}
\label{sep-non-sep-cicrles}
\end{figure}

\myitpar{Two examples.} The system of circles\dss $s$\dss pictured on Fig.\dss 
\ref{separating-circles}\dss consists
of\qss $2\dff\gen\qff -\qff 3$\qss separating circles.\oss
Therefore,\oss for this system of circles\qss $N\qff =\qff 2\dff\gen\qff -\qff 3$\snsp,\oss
$n\qff =\qff 0$\snsp,\oss and hence\oss
$\displaystyle
\rank\dff \mathcal{T}\fff(s)\dff\cap\dff \tors\off =\off 2\dff\gen\qff -\qff 3$\snsp.\oss
The system of circles\dss $s'$\dss pictured on Fig.\dss \ref{sep-non-sep-cicrles}\dss consists
of\qss $\gen\qff -\qff 1$\qss separating circles and\qss 
$\gen\qff -\qff 1$\qss non-separating circles,\oss all of which are homology equivalent\halfff.\oss
Therefore,\oss for this system of circles\qss $N\qff =\qff 2\dff\gen\qff -\qff 2$\snsp,\oss
$n\qff =\qff 1$\snsp,\oss and hence\oss
$\displaystyle
\rank\dff \mathcal{T}\fff(s')\dff\cap\dff \tors\off =\off 
\left(2\dff\gen\qff -\qff 2\right)\qff -\qff 1\off =\off
2\dff\gen\qff -\qff 3$\snsp.\oss

\renewcommand{\topfraction}{1.0}
\renewcommand{\textfraction}{0.0}

\renewcommand{\topfraction}{1.0}
\renewcommand{\textfraction}{0.0}

\mypar{Lemma.}{number-of-necklaces} \emph{Let\qss $s$\qss 
be a one-dimensional closed submanifold of\qss $S$\qss
partitioning\qss $S$\qss into parts of genus\dss $0$\dnsp.\oss 
Then\dss $s$\dss has at least\pss $\genus(S)$\dss necklaces.\oss}

\proof\qss Let argue by induction by the genus\qss $\gen\qff =\qff \genus(S)$\snsp.\oss
The cases of\qss $\gen\qff =\qff 0\fff,\pff 1$\qss are trivial.\oss
Let us consider a component\qss $C$\qss of\qss $s$\dnsp.\oss
Let us cut\qss $S$\qss along\qss $C$\qss and glue 
two discs to the two boundary components of the resulting surface
and denote by\qss $R$\qss the resulting surface.\oss
The submanifold\qss $s\smallsetminus C$\qss partitions\qss $R$\qss into discs with holes.\oss

If\qss $C$\qss is separating\halfff,\oss then\qss $R$\qss consists of two components.\qss
Let us denote them by\qss $R_{\fff 1}$\qss and\qss $R_{\fff 2}$\nsp,\oss
and\dss let\qss $\gen_{\fff 1}\qff =\qff \genus(R_{\fff 1})$\qss and\qss
$\gen_{\fff 2}\qff =\qff \genus(R_{\fff 2})$\dnsp.\oss
Then\qss $\gen\qff =\qff \gen_{\fff 1}\qff +\qff \gen_{\fff 2}$\qss
and\qss $\gen_{\fff 1}\fff,\pff \gen_{\fff 2}\dff <\dff \gen$\snsp.\oss
The submanifolds\oss
$\displaystyle s_1\qff =\qff (s\smallsetminus C)\dff\cap\dff R_{\fff 1}$\oss and\oss
$\displaystyle s_2\qff =\qff (s\smallsetminus C)\dff\cap\dff R_{\fff 2}$\oss
partition\qss $R_{\fff 1}$\qss and\qss $R_{\fff 2}$\qss respectively into parts of genus\dss $0$.\oss
By the inductive assumption submanifolds\dss $s_1$\dss and\dss $s_2$\dss have  
at least\dss $\gen_{\fff 1}$\dss and\dss $\gen_{\fff 2}$\dss necklaces respectively.\oss
Each of these necklaces is also a necklace of\dss $s$\dss in\qss $S$\dnsp.\qff\oss
It follows that\dss $s$\dss has at least\qss $\gen$\qss necklaces.\oss

If\qss $C$\qss is non-separating\halfff,\oss then\qss $R$\qss is connected and\oss
$\displaystyle
\genus(R)\qff =\qff \gen\qff -\qff 1$\dnsp.\oss
Since\qss $C$\qss is non-separating\halfff,\qss $C$\qss belongs to some necklace of\qss $s$\dnsp.\oss
Other components of this necklace are separating in\dss $R$\dss
and hence do not belong to any necklace of\qss $s\smallsetminus C$\qss in\dss $R$\snsp.\oss
On the other hand,\oss any non-separating circle in\dss $R$\qss 
is non-separating in\dss $S$\snsp,\pss and if two such circles are not homology equivalent
in\dss $R$\nsp,\oss then they are not homology equivalent in\qss $S$\dnsp.\oss 
Hence the number of necklaces of\qss $s\smallsetminus C$\qss in\qss $R$\qss
is smaller 
than the number of necklaces of\dss $s$\dnsp.\oss
By the inductive assumption there are at least\qss $\gen\qff -\qff 1$\qss necklaces of\qss
$s\smallsetminus C$\qss in\qss $R$\snsp,\oss
and hence at least\qss $\gen$\qss necklaces of\trs $s$\dss in\qss $S$\dnsp.\oss
This completes the step of the induction.\oss  \eproof

\mypar{Corollary.}{multi-twists-subgroups-general} \emph{Let\dss $s$\dss 
be a one-dimensional closed submanifold of\qss $S$\dnsp.\qff\oss
Then} 
\[
\quad
\rank\dff \mathcal{T}\fff(s)\dff\cap\dff \tors\off \leq\off 2\dff\gen\qff -\qff 3\fff.
\]

\vspace*{-\bigskipamount}
\proof\qss Since Dehn twists about trivial circles are equal to $1$
and Dehn twists about isotopic circles are equal,\oss
we may assume that\dss $s$\dss is a system of circles.\qff\oss
Adding new components to\dss $s$\dss cannot decrease the rank of\qss
$\mathcal{T}\fff(s)\dff\cap\dff \tors$\dnsp.\oss
Therefore,\oss we may assume that\dss $s$\dss partitions\qss $S$\qss into discs with two holes.\oss
Then the number of components of\dss $s$\dss is equal to\qss $3\dff\gen\qff -\qff 3$\dnsp,\oss
and\dss Corollary\qss \ref{multi-twists-subgroups}\qss
together with Lemma\qss \ref{number-of-necklaces}\qss imply that
the rank of\oss $\displaystyle \mathcal{T}\fff(s)\dff\cap\dff \tors$\oss is\oss
$\leq\off 3\dff\gen\qff -\qff 3\qff -\qff \gen\off =\off 2\dff\gen\qff -\qff 3$\dnsp.\oss  \eproof

\myitpar{An invariant of surfaces with boundary.} Let\qss $Q$\qss
be a compact orientable surface which is not an annulus
and which has non-empty boundary.\oss
As usual,\oss let\qss $\gen\qff =\qff \genus(Q)$\qss be the genus of\qss $Q$\dnsp.\oss
Let\qss $b\qff =\qff b\fff(Q)\qff \geq\qff 1$\qss 
be the number of components of\qss $\partial Q$\dnsp.\oss
Let\oss
\[
\quad
d\fff(Q)\off =\off 
2\dff\genus(Q)\qff -\qff 3\qff +\qff b\fff(Q)
\off =\off 
2\dff\gen\qff -\qff 3\qff +\qff b\dff.
\]
The geometric meaning of\dss $d\fff(Q)$\dss is the following.\oss
The maximal number of components of a system of circles\dss $c$\dss on\qss $Q$\qss
is equal to\qss $3\dff\gen\qff -\qff 3\qff +\qff b$\dnsp,\oss
and any system of circles with the maximal number of components
partitions\qss $Q$\qss into discs with two holes.\oss
One can construct system of circles with maximal number of components,\oss
in particular\halfff,\oss as follows.\oss
Choose first\qss $\gen$\qss disjoint non-trivial circles in\qss $Q$\qss such that 
their union\dss $c_{\fff 0}$\dss does separate\qss $Q$\dnsp.\oss
The surface\qss $Q$\qss cut along\qss $c_{\fff 0}$\qss has genus\dss $0$\dnsp.\oss
By adding\dss $d\fff(Q)$\dss circles to\dss $c_{\fff 0}$\dss
one can get a system of circles on\qss $Q$\qss with the maximal number of components.\oss

\myitpar{Two invariants of systems of circles.} Let\dss $s$\dss
be a system of circles on\trs $S$\dnsp.\off\oss 
Let\oss\vspace*{-0.5\medskipamount}
\[
\quad
\mathfrak{D}(s)\off\dff =\off\qff \sum\nolimits_{\dff Q}\qff d\fff(Q)\dff,
\] 

\vspace*{-1.25\bigskipamount}
where\qss $Q$\qss runs over components of\pss $S\ccut s$\snsp.\qff\oss 
Since\dss $s$\dss is a system of circles,\oss
none of components\qss $Q$\qss of\qss $S\ccut s$\qss is an annulus,\oss
and hence\qss $d\fff(Q)$\qss is defined for all\qss $Q$\dnsp.\oss

Let\qss $\mathfrak{d}(s)$\qss be the number of components of\pss 
$S\ccut s$\pss which are neither a disc with two holes,\oss
nor a torus with one hole.\oss
Then\oss\vspace*{-0.5\medskipamount}
\[
\quad
\mathfrak{D}(s)\off \geq\off \mathfrak{d}(s)\fff
\]

\vspace*{-1.25\bigskipamount}
because\qss $d\fff(Q)\qff \geq\qff 1$\qss if\qss $Q$\qss
is neither a disc with two holes,\oss nor a torus with one hole.\oss

\mypar{Corollary.}{stronger-general-bound} \emph{Let\dss $s$\dss 
be a system of circles on\trs $S$\dnsp.\off\oss 
Then}
\[
\quad
\rank\dff \mathcal{T}\fff(s)\dff\cap\dff \tors
\off \leq\off
2\dff\gen\qff -\qff 3\qff -\qff \mathfrak{D}(s)
\off \leq\off 
2\dff\gen\qff -\qff 3\qff -\qff \mathfrak{d}(s)\fff.
\]

\vspace*{-\bigskipamount}
\proof\qss Each component\dss $Q$\dss of\qss $S\ccut s$\qss with\qss $\genus(Q)\qff \geq\qff 1$\qss
contains\qss $\genus(Q)$\qss disjoint circles such that 
their union does separate\qss $Q$\dnsp.\oss
All these circles are pair-wise not homology equivalent\halfff,\oss 
and not homology equivalent to any component of\dss $s$\dnsp.\oss
Let\qss $s_{\fff 1}$\qss be the union of\qss $s$\qss and all these circles.\oss
Then\qss $s_{\fff 1}$\qss is a system of circles,\oss
and every new circle in\qss $s_{\fff 1}$\qss 
is the single element of a new necklace of\qss $s_{\fff 1}$\dnsp.\oss
Therefore Theorem\qss \ref{multi-twists-in-t}\qss implies that\oss
\[
\quad
\mathcal{T}\fff(s_{\fff 1})\dff\cap\dff \tors\off 
=\off \mathcal{T}\fff(s)\dff\cap\dff \tors\dff.
\]
Let\qss $N_{\fff 1}$\qss is the number of components of\pss $s_{\fff 1}$\snsp,\oss
and\qss $n_{\fff 1}$\qss is the number of necklaces of\pss $s_{\fff 1}$\snsp.\oss
By the construction,\qss $s_{\fff 1}$\qss partitions\qss $S$\qss 
into subsurfaces of genus\dss $0$\dnsp.\oss
Hence Lemma\qss \ref{number-of-necklaces}\qss implies
that there are at least\dss $\gen$\dss 
necklaces of\qss $s_{\fff 1}$\dnsp,\qff\oss
i.e.\qss $n_{\fff 1}\qff \geq\qff \gen$\nsp.\oss
One can get from\qss $s_{\fff 1}$\qss a system of circles in\qss $S$\qss
with the maximal possible number of components
by adding\dss $d\fff(Q)$\dss circles in\qss $Q$\qss for each component\qss $Q$\dnsp.\oss
Since the maximal possible number of components is\qss $3\dff\gen\qff -\qff 3$\dnsp,\oss
it follows that\qss $N_{\fff 1}\qff +\qff \mathfrak{D}(s)\qff \leq\qff 3\dff\gen\qff -\qff 3$\dnsp,\oss
and hence\qss $N_{\fff 1}\qff \leq\qff 3\dff\gen\qff -\qff 3\qff -\qff \mathfrak{D}(s)$\dnsp.\oss
By combining the last inequality with\qss $n_{\fff 1}\qff \geq\qff \gen$\qss 
we see that
\[
\quad
N_{\fff 1}\qff -\qff n_{\fff 1}
\off \leq \off
3\dff\gen\qff -\qff 3\qff -\qff \mathfrak{D}(s)\qff -\qff \gen
\off =\off
2\dff\gen\qff -\qff 3\qff -\qff \mathfrak{D}(s)\dff.
\]
But by Corollary\qss \ref{multi-twists-subgroups}\oss
$\displaystyle
\mathcal{T}\fff(s_{\fff 1})\dff\cap\dff \tors$\oss is a free abelian group of rank\qss
$N_{\fff 1}\qff -\qff n_{\fff 1}$\dnsp.\oss  \eproof

\mysection{Pure\qss diffeomorphisms\qss and\qss reduction\qss systems}{pure-reduction}

\myitpar{Pure diffeomorphisms and elements.} 
Let\qss $F$\qss be a diffeomorphism  of\qss $S$\dnsp.\oss
A system of circles\dss $c$\dss is said to be a\qss
\emph{pure reduction system}\qss for\dss $F$\qss if\dss $c$\dss is a reduction system for\qss $F$\qss
and the following four conditions hold.\vspace*{-\medskipamount}
\begin{itemize}
\item[(a)] $F$\qss is orientation-preserving\fff.\oss 
\item[({\fff}b)] {\qff}Every component of\qss $S\cutc$\qss is invariant under\trs $F\cutc$\snsp.\oss  
\item[(c)] $F$\qss is equal to the identity 
           in a neighborhood of\dss $c\dff\cup\dff\partial S$\snsp.\oss
\item[(d)] {\qff}For each component\trs $Q$\trs of\qss $S\cutc$\qss the isotopy class of 
the restriction\qss $F_Q\fff\colon Q\qff\toto\qff Q$\qss
{\qff}is either pseudo-Anosov,\oss or\dss contains\qss $\id_Q$\snsp.\oss
\end{itemize}

\vspace*{-\medskipamount}
A diffeomorphism\dss $F$\dss of\dss $S$\dss is said to be\pss \emph{pure}\oss 
if\qss $F$\qss admits a pure reduction system.\oss
This definition is invariant under diffeomorphisms of\qss $S$\dnsp.\oss
More formally,\oss if\dss $c$\dss is a pure reduction system 
for a diffeomorphism\dss $F$\dss of\dss $S$\dnsp,\oss
and if\qss $G$\qss is some other diffeomorphism of\qss $S$\dnsp,\oss
then\qss $F\dff(c)$\dss is a pure reduction system for\qss $G\circ F\circ G^{\fff -\halfff 1}$\snsp.\oss

An isotopy class\dss $f\dff\in\dff\mmod(S)$\dss is\dss said\dss to\dss be\pss \emph{pure}\fff\oss 
if\fff\qss $f$\qss contains\dss a\dss pure\dss diffeomorphism.\oss
A system of circles\dss $c$\dss is said to be a\qss
\emph{pure reduction system}\qss for an element\qss $f\dff\in\dff \mms$\qss
if\dss $c$\dss is a pure reduction system for some diffeomorphism\qss $F$\qss
in the isotopy class\qss $f$\snsp.\oss

The\qss \emph{isotopy extension theorem}\qss implies that the property of being a pure
reduction system for\qss $f$\qss depends only 
on the isotopy class of the submanifold\dss $c$\dss in\qss $S$\dnsp.\oss
In addition,\oss if\dss $c$\dss is a pure reduction system for\dss $f$\dss
and\dss $g$\dss is the isotopy class of a diffeomorphism\dss $G$\dss of\dss $S$\dnsp,\oss 
then\qss $G\fff(c)$\qss is a pure reduction system for\qss $g\dff f\dff g^{\fff -\dff 1}$\dnsp.\oss

\mypar{Theorem.}{pure-torelli} \emph{If\pss $m\qff\geq\qff 3$\snsp,\qff\oss
then all elements of\pss $\ttt_{\fff m\fff}(S)$\pss are pure.}

\mypar{Lemma.}{pure-abc} \emph{Suppose that\dss $c$\dss 
is a reduction system for a diffeomorphism\qss $G$\qss
representing an element of\qss $\ttt_{\fff m\fff}(S)$\dnsp,\oss 
where\pss $m\qff\geq\qff 3$\snsp.\qff\oss
Then\qss $G$\pss leaves every component of\qss $c$\qss invariant\halfff,\oss
and\qss $G\cutc$\pss leaves every component of\oss $S\cutc$\qss invariant\halfff.\oss}

\prooftitle{Proofs}\qss See\qss Theorem\qss 1.7\qss and\qss Theorem\qss 1.2\qss 
of\oss \cite{i-book}\oss respectively.\oss  \eproof

\mypar{Lemma.}{boundary-pseudo-anosov} \emph{Suppose that\pss 
$f\dff\in\dff \mms$\qss is a pure element\halfff.\oss
Let\dss $c$\trs be a pure reduction system for a diffeomorphism\qss $F\fff\colon S\ttoo S$\qss
representing\qss $f$\snsp.\qff\oss
If\oss $Q$\oss is a component of\oss $S\cutc$\pss such that 
the restriction\pss $F_Q$\qss is pseudo-Anosov,\oss
then the image\pss $p\ccut c\dff(\partial Q)\dff\subset\dff c$\pss is contained up to isotopy in any
pure reduction system of\pss $f$\snsp.\oss}

\proof\qss
This follows from the uniqueness of Thurston's normal form of\qss $f$\snsp.\oss  \eproof

\myitpar{Minimal pure reduction systems.} A system of circles\dss $c$\dss is said to be a\qss
\emph{minimal pure reduction system}\qss for an element\qss $f\dff\in\dff \mms$\qss if\dss
$c$\dss is a pure reduction system for\qss $f$\snsp,\oss but no proper subsystem of\dss $c$\dss is.\oss
Any pure reduction system for\qss $f$\qss contains a minimal pure reduction system,\oss
and,\oss in fact\halfff,\oss it is unique.\oss

Moreover\halfff,\oss up to isotopy such a minimal pure reduction system 
depends only on\qss $f$\snsp.\oss

In fact\halfff,\oss the set of the isotopy classes of 
components of a minimal pure reduction system for\qss $f$\qss
is nothing else but the canonical reduction system of\qss $f$\qss 
in the sense of\qss \cite{i-book}.\oss 
This easily follows from the results of\oss \cite{i-book}\halfff,\oss Chapter 7.\oss
See also\oss \cite{im},\oss Section\qss 3.\oss
Since the canonical reduction system of\qss $f$\qss 
is defined invariantly in terms of\qss $f$\snsp,\oss
it depends only on\qss $f$\qss and hence the same is true 
for the minimal reduction systems for\qss $f$\snsp.\oss

Moreover\halfff,\oss the notion of a minimal pure reduction system is invariant
under diffeomorphisms of\qss $S$\qss in the same sense as the notion of
a pure reduction system\qss (see the first subsection of this section).\oss

\myitpar{Reduction systems of subgroups.} Let\dss $\Gamma$\dss be a subgroup of\qss $\mms$\dnsp.\oss
A system of circles\dss $c$\dss on\dss $S$\dss 
is called a\qss \emph{reduction system}\qss 
for\dss $\Gamma$\dss if\dss $c$\dss is a reduction system 
for every\qss $f\dff\in\dff \Gamma$\dnsp.\oss 
The subgroup\dss $\Gamma$\dss is said to be\qss \emph{reducible}\qss 
if\dss there exists a non-empty reduction system for\qss $\Gamma$\dnsp,\oss
and\dss is\sss said to be\qss \emph{irreducible}\qss otherwise.\oss

A finite subgroup of\qss $\mmod(S)$\qss can be irreducible even if all its
non-trivial elements are reducible.\qss
Such subgroups were constructed and classified by\dss J.\dss Gilman\qss \cite{gilman}.\oss
But an infinite irreducible subgroup always contains an irreducible element
of infinite order{\halfff},\oss i.e.\trs a pseudo-Anosov element.\oss
See\oss \cite{i-book},\oss Corollary\oss 7.14.\oss

\myitpar{Pure reduction systems of subgroups.} 
Suppose that\qss $\Gamma$\qss is a subgroup of\qss $\mms$\qss
consisting of pure elements.\oss 
A system of circles\dss $c$\dss is said to be a\qss
\emph{pure reduction system}\qss for\qss $\Gamma$\qss if\dss $c$\dss
is a pure reduction system for every element of\qss $\Gamma$\snsp.\oss

In general,\oss a subgroup of\qss $\mms$\qss consisting of pure elements
does not admit any pure reduction system.\qss
For example,\oss $\tors$\qss does not admit a pure reduction system.\oss
Indeed,\oss $\tors$\qss contains both pseudo-Anosov elements and Dehn multi-twists.\oss
On the other hand,\oss a reduction system of a pseudo-Anosov element should be empty,\oss
but a reduction system of a Dehn multi-twist cannot be empty.\oss

\mypar{Theorem.}{abelian-reduction} \emph{If\pss $\mathfrak{G}$\pss 
is an abelian subgroup of\pss $\ttt_{\fff m\fff}(S)$\dnsp,\oss 
where\qss $m\qff\geq\qff 3$\snsp,\qff\oss
then there exists a pure reduction system for\qss $\mathfrak{G}$\dnsp.\oss
In particular\halfff,\oss if\pss $\mathfrak{G}$\pss 
is an abelian subgroup of\pss $\tors$\dnsp,\qff\oss
then there exists a pure reduction system for\qss $\mathfrak{G}$\dnsp.\oss}

\proof\dss See the description of abelian subgroups of\dss $\mms$\dss
in\qss \cite{i-book}\halfff,\oss Section\qss 8.12.\oss  \eproof

\myitpar{Abelian subgroups of\qss $\tors$\dnsp.} Let\qss $\mathfrak{A}$\qss 
be an abelian subgroup of the group\pss 
$\ttt_{\fff m\fff}(S)$\dnsp,\oss where\pss $m\qff\geq\qff 3$\snsp.\qff\oss
In particular\halfff,\oss $\mathfrak{A}$\qss may be 
an abelian subgroup of\pss $\tors$\dnsp.\qff\oss
Let\dss $c$\dss be a pure reduction system for\qss $\mathfrak{A}$\nsp.\oss
Then\trs $\mathfrak{A}$\trs is contained 
in a free abelian subgroup\trs $\mathfrak{G}$\trs of\qss $\mms$\qss
constructed as follows.

Suppose that\dss $\pi$\dss is a set of components of\qss $S\cutc$\dnsp.\oss 
Suppose that for each component\dss $Q\dff\in\dff \pi$\dss
a diffeomorphisms\qss $F^{\fff Q}\dff\colon Q\ttoo Q$\qss fixed on\dss $\partial Q$\qss is given.\qff\oss 
Let\qss $f^{\dff Q}\dff\in\dff \mm(Q)$\qss be the isotopy class of\qss$F^{\fff Q}$\nsp.\oss
Suppose that each isotopy class\dss $f^{\dff Q}$\dss is pseudo-Anosov.\oss
Let us extend these diffeomorphisms by the identity to diffeomorphisms of\dss $S$\dnsp,\oss
and let\dss $\mathfrak{G}$\dss be the subgroup of\qss $\mms$\qss generated by\qss
$\mathcal{T}\fff(s)$\qss and the isotopy classes of these extensions.\oss

Then\dss $\mathfrak{G}$\dss is an abelian subgroup of\qss $\mms$\dnsp,\oss 
and\pss $\mathfrak{A}\dff\subset\dff \mathfrak{G}$\dnsp.\oss
Moreover\halfff,\oss these extensions and the Dehn twists about components of\dss $c$\dss
are free generators of\qss $\mathfrak{G}$\dnsp.\oss

This is just a rephrased description of abelian subgroups
from\qss \cite{i-book}\halfff,\oss Section\qss 8.12.\oss

\mysection{The\qss rank\qss of\qss abelian\qss subgroups\qss of\qss Torelli\qss groups}{abelian-torelli}

\mypar{Theorem.}{abelian-rank} \emph{Every abelian subgroup 
of\pss $\tors$\qss is a free abelian group of rank\oss 
$\leq\off 2\dff\gen\qff -\qff 3$\dnsp.\oss
If\dss $c$\dss is a pure reduction system 
for an abelian subgroup\oss $\mathfrak{A}\dff\subset\dff \tors$\dnsp,\qff\oss
then\oss}
\[
\quad
\rank \mathfrak{A}
\off \leq\off 
\mathfrak{d}(c)\qff +\qff \rank\qff \mathcal{T}\fff(c)\fff\cap\dff\tors
\off \leq\off
2\dff\gen\qff -\qff 3\qff -\qff \bigl(\dff\mathfrak{D}(c)\qff -\qff \mathfrak{d}(c)\dff\bigr)\fff.
\]

\vspace*{-\bigskipamount}
\proof\qss Let\dss $c$\dss be a pure reduction system for 
an abelian subgroup\trs $\mathfrak{A}$\trs of\qss $\tors$\dnsp.\oss
The subgroup\trs $\mathfrak{A}$\trs is contained in a free abelian subgroup\qss
$\mathfrak{G}$\qss of\qss $\mms$\qss 
of the form described at the end of\dss Section\qss \ref{pure-reduction}.\oss
In the rest of the proof we will use notations introduced in the construction of\trs
$\mathfrak{G}$\trs in\dss Section\qss \ref{pure-reduction}.\oss

Restriction to the components of\qss $S\cutc$\qss defines a canonical surjective homomorphism
\[
\quad
\rho\dff\colon\dff \mathfrak{G}\ttoo \prod\nolimits_{\dff Q}\qff \zzz^{\fff Q}\fff,
\]
where the product is taken over components\dss $Q\dff\in\dff \pi$\nsp,\oss 
and\qss $\zzz^{\fff Q}$\qss is the infinite cyclic 
subgroup of\qss $\mm(Q)$\qss generated by\qss $f^{\dff Q}$\dnsp.\oss
The image of\dss $\rho$\dss is a free abelian group of rank equal to the number
of elements of\dss $\pi$\nsp,\oss 
and the kernel of\qss $\rho$\qss is equal to\qss $\mathcal{T}\fff(s)$\dnsp.\qff\oss
Let\oss 
$\rho\fff|\dff_\mathfrak{A}$\oss
be the restriction of\qss $\rho$\qss to\qss $\mathfrak{A}$\dnsp.\oss
Let us estimate the ranks of its image and kernel.\oss

If\qss $Q\dff\in\dff \pi$\dnsp,\oss then\qss $Q$\qss is not a disc with two holes
because\qss $f^{\dff Q}$\qss is a pseudo-Anosov class.\oss
If\qss $Q\dff\in\dff \pi$\qss and\qss $Q$\qss is a torus with one hole,\oss
then\qss $F^{\fff Q}$\qss acts non-trivially on\qss $H_{\fff 1}(Q)$\dnsp.\oss
In this case\qss $\partial Q$\qss is a separating circle in\qss $S$\dnsp,\oss
and hence\qss $H_{\fff 1}(Q)$\qss is a direct summand of\qss $H_{\fff 1}(S)$\dnsp.\oss
It follows that\qss $F^{\fff Q}$\qss cannot be the restriction to\qss $Q$\qss
of a diffeomorphism\qss $S\ttoo S$\qss acting trivially on\qss $H_{\fff 1}(S)$\dnsp.\oss
In turn,\oss this implies that the image of\qss $\rho\fff|\dff_\mathfrak{A}$\qss
is contained in the product of factors\qss $\zzz^{\fff Q}$\qss with\qss $Q$\qss 
not a torus with one hole.\oss 

It follows that the rank of the image of\qss $\rho\fff|\dff_\mathfrak{A}$\qss 
is\oss 
$\displaystyle \leq\off \mathfrak{d}(c)$\dnsp.\oss 
The kernel of\qss $\rho\fff|\dff_\mathfrak{A}$\qss  is contained in\qss 
$\displaystyle
\mathcal{T}\fff(s)\fff\cap\dff\mathfrak{A}\off\subset\off \mathcal{T}\fff(s)\fff\cap\dff\tors$\dnsp.\oss
It follows that
\[
\quad
\rank \mathfrak{A}
\off \leq\off 
\mathfrak{d}(c)\qff +\qff \rank\qff \mathcal{T}\fff(c)\fff\cap\dff\tors\fff.
\]
By Corollary\qss \ref{stronger-general-bound}\oss 
$\displaystyle
\rank\qff \mathcal{T}\fff(s)\fff\cap\dff\tors
\off \leq\off 2\dff\gen\qff -\qff 3\qff -\qff \mathfrak{D}(c)$\nsp,\oss 
and hence\oss 
\[
\quad
\rank \mathfrak{A}
\off \leq\off 
\mathfrak{d}(c)\qff +\qff \bigl(\dff2\dff\gen\qff -\qff 3\qff -\qff \mathfrak{D}(c)\dff\bigr)
\off =\off
2\dff\gen\qff -\qff 3\qff -\qff \bigl(\dff\mathfrak{D}(c)\qff -\qff \mathfrak{d}(c)\dff\bigr)\fff.
\]
Since\oss $\mathfrak{D}(c)\qff \geq\qff \mathfrak{d}(c)$\dnsp,\oss
this implies that\qss $\rank \mathfrak{A}\qff \leq\qff 3\dff\gen\qff -\qff 3$\dnsp.\oss  \eproof

\mypar{Lemma.}{number-of-necklaces-better} \emph{Let\dss $c$\dss 
be a system of circles
partitioning\qss $S$\qss into surfaces of genus\dss $0$\dnsp.\qff\oss 
If\dss there\sss is\sss a\sss component\dss $Q$\dss of\pss $S\cutc$\pss
such that\trs $\genus(Q)\qff =\qff 0$\dnsp,\oss 
the canonical map\qss
$p\cutc$\qss embeds\trs $Q$\trs in\qss $S$\dnsp,\oss 
and the complementary surface\qss $\ccomp Q$\qss is connected,\oss
then\dss $c$\dss has at least\pss $\genus(S)\qff +\qff 1$\pss necklaces.\oss}

\proof\qss Let us argue by induction by the number of components of\qss $\partial Q$\dnsp.\oss
Since\dss $c$\dss is a system of circles,\oss $Q$\qss is neither a disc,\oss
nor an annulus,\oss and hence this number is\qss $\geq 3$\dnsp.\oss

Suppose that\qss $\partial Q$\qss consists of\dss $3$\dss components,\oss
and let\qss $C\fff,\pff D\fff,\pff E$\qss be these components.\oss 
Since\qss $\ccomp Q$\qss is connected,\oss all of them are non-separating
and no two of them are homology equivalent\halfff.\oss
Hence\qss $C\fff,\pff D\fff,\pff E$\qss belong to three different necklaces.\oss

As in the proof of Lemma\qss \ref{number-of-necklaces},\oss
let us cut\qss $S$\qss along\qss $C$\qss and glue 
two discs to the two boundary components of the resulting surface.\oss
Let\qss $R$\qss be the result of this glueing.\oss
The necklace containing\qss $C$\qss disappears in\qss $R$\qss
(since the two circles in\qss $R$\qss resulting from\qss $C$\qss bound discs in\qss $R$\nsp),\oss
and the two necklaces containing\qss $D$\qss and\qss $E$\qss respectively
coalesce in\qss $R$\qss into one\oss (since the union\dss $D\dff\cup\dff E$\qss
bounds in\qss $R$\qss an annulus).\oss
Hence the number of necklaces of\dss $c$\dss in\qss $S$\qss 
is bigger than the number of necklaces of\qss 
$c\smallsetminus C$\qss in\qss $R$\qss by at least\dss $2$\dnsp.\oss
By Lemma\qss \ref{number-of-necklaces}\qss there are at least\qss
$\genus(R)\qff =\qff\dff \genus(S)\qff -\qff 1$\qss necklaces of\dss 
$c\smallsetminus C$\qss in\qss $R$\dnsp,\oss
and hence there are at least\oss
$\genus(S)\qff +\qff 1$\oss
necklaces\dss of\fff\qss $c$\dss in\qss $S$\dnsp.\oss  

Suppose now that\trs $\partial Q$\trs consists of\qss $\geq\qff 4$\qss components,\oss
and let\qss $C$\qss be one of them.\oss
Again,\oss let us cut\qss $S$\qss along\qss $C$\qss and glue 
two discs to the two boundary components of the resulting surface,\oss
and let\qss $R$\qss be the result of this glueing\halfff.\oss
One of the two glued discs is glued to\trs $Q$\dnsp.\oss
Let\qss $P$\qss be the result of this glueing\halfff.\qff\oss
Let\qss $c'\qff =\qff c\smallsetminus C$\dnsp.\qff\oss
Then\trs $c'$\trs is a system of circles on\qss $R$\snsp,\oss
and\dss $P$\dss is a component of\dss $R\dff\ccut c'$\dnsp.\oss
Moreover\halfff,\oss $p\ccut c'$\qss embeds\qss $P$\qss in\qss $R$\snsp,\oss
and the complementary surface\qss $\ccomp P$\qss is connected\qss
({\fff}being the result of glueing a disc to a boundary component of\qss $\ccomp Q$\nsp).\oss

The necklace\dss $\mathfrak{n}_{\dff C}$\dss containing\dss $C$\dss disappears in\dss $R$\dss
({\fff}by the same reason as above).\oss
If a non-separating component of\dss $c$\dss does not belong to\dss $\mathfrak{n}_{\dff C}$\nsp,\oss
then it is non-separating in\dss $R$\dss also,\oss
and every non-separating in\dss $R$\dss component 
of\trs $c'$\trs is non-separating in\qss $S$\qss also.\oss
If two non-separating components of\dss $c$\dss do not belong to\dss $\mathfrak{n}_{\dff C}$\dss
and are homology equivalent in\qss $S$\nsp,\oss 
then they are homology equivalent in\qss $R$\qss also.\oss 

It follows that there is a canonical surjective map from the set of  
different from\dss $\mathfrak{n}_{\dff C}$\dss necklaces of\dss $c$\dss in\qss $S$\qss
to the set of necklaces of\trs $c'$\trs in\qss $R$\nsp.\oss
Hence the number of necklaces of\dss $c$\dss in\qss $S$\qss 
is bigger than the number of necklaces of\trs 
$c'$\trs in\qss $R$\qss by at least\dss $1$\dnsp.\oss
By the inductive assumption there are at least\qss 
$\genus(R)\qff +\qff 1\qff =\qff\dff \genus(S)$\qss
necklaces of\trs $c'$\snsp.\oss
It follows that there aree at least\qss 
$\genus(S)\qff +\qff 1$\qss necklaces of\trs $c$\snsp.\dff\oss
This completes the induction step and hence the proof of the lemma.\oss  \eproof

\mypar{Theorem.}{lower-rank-necklaces} \emph{Let\dss $c$\dss be a pure reduction system 
of an abelian subgroup\pss $\mathfrak{A}\dff\subset\dff \tors$\dnsp.\qff\oss
Suppose that there\sss is\sss a\sss component\qss $Q$\qss of\oss $S\cutc$\pss
such that\pss $\genus(Q)\qff =\qff 0$\dnsp,\qff\oss 
the canonical map\qss $p\cutc$\qss embeds\trs $Q$\qss in\qss $S$\dnsp,\oss 
and the complementary surface\qss $\ccomp Q$\qss is connected.\qff\oss
Then}\oss $\rank \mathfrak{A}\off \leq \off 2\dff\gen\qff -\qff 4$\dnsp.\oss

\proof\qss If a component\dss $Q$\dss of\qss $S\cutc$\qss is neither a disc with two holes,\oss 
nor a torus with one hole,\oss
then there is a circle contained in\dss $Q$\dss and non-trivial in\dss $Q$\dnsp.\oss
By adding these circles to\dss $c$\dss we will get a new system of circles\trs $c'$\snsp.\oss
By the definition of\qss $\mathfrak{d}(c)$\dnsp,\oss 
there are\qss $\mathfrak{d}(c)$\qss such components,\oss
and hence there are\qss $\mathfrak{d}(c)$\qss new circles in\trs $c'$\snsp.\oss
On the other hand,\oss $c'$\trs consists of\qss $\leq\qff 3\dff\gen\qff -\qff 3$\qss
components because\trs $c'$\trs is a system of circles.\oss
It follows that the number of components of\dss $c$\dss is\qss
$\leq\qff 3\dff\gen\qff -\qff 3\qff -\qff \mathfrak{d}(c)$\dnsp.\oss
On the other hand,\oss
Lemma\qss \ref{number-of-necklaces-better}\qss implies that under our assumptions
the number of necklaces of\qss $c$\qss is\qss $\geq\qff g\qff +\qff 1$\dnsp.\oss

By combining these estimates of the number of components and the number of necklaces of\dss $c$\dss
with\dss Theorem\qss \ref{multi-twists-subgroups}\fff,\oss we see that\oss
\[
\quad
\rank\qff \mathcal{T}\fff(c)\fff\cap\dff\tors
\off \leq\off
3\dff\gen\qff -\qff 3\qff -\qff \mathfrak{d}(c)\qff -\qff (\gen\qff +\qff 1)
\off =\off
2\dff\gen\qff -\qff 4\qff -\qff \mathfrak{d}(c)\fff.
\]
By\dss Theorem\qss \ref{abelian-rank}\oss 
$\displaystyle
\rank \mathfrak{A}
\off \leq\off 
\mathfrak{d}(c)\qff +\qff \rank\qff \mathcal{T}\fff(c)\fff\cap\dff\tors$\dnsp,\qff\oss
and\dss hence\oss
\[
\quad
\rank \mathfrak{A}
\off \leq\off 
\mathfrak{d}(c)\qff +\qff \bigl(\dff 2\dff\gen\qff -\qff 4\qff -\qff \mathfrak{d}(c)\dff\bigr)
\off =\off
2\dff\gen\qff -\qff 4\fff.
\]
This completes the proof of the theorem.\oss  \eproof

\mypar{Theorem.}{lower-rank-components} \emph{Let\dss $c$\dss be a pure reduction system 
of an abelian subgroup\oss $\mathfrak{A}\dff\subset\dff \tors$\dnsp.\qff\oss
Then}\oss $\rank \mathfrak{A}\off \leq \off 2\dff\gen\qff -\qff 4$\oss
\emph{unless each component of\oss $S\cutc$\pss is either a sphere with\trs $3$\trs or\trs $4$\trs holes,\oss
or a torus with\trs $1$\trs or\trs $2$\qss holes.\oss}

\proof\qss By Theorem\qss \ref{abelian-rank}\oss
$\displaystyle
\rank \mathfrak{A}
\off \leq\off
2\dff\gen\qff -\qff 3\qff -\qff \bigl(\dff\mathfrak{D}(c)\qff -\qff \mathfrak{d}(c)\dff\bigr)$\dnsp.\oss
Since\qss $\mathfrak{D}(c)\qff \geq\qff \mathfrak{d}(c)$\dnsp,\oss
\[
\quad
\rank \mathfrak{A}
\off \leq\off
2\dff\gen\qff -\qff 4
\]
unless\qss $\mathfrak{D}(c)\qff =\qff \mathfrak{d}(c)$\dnsp.\oss
If the last equality holds,\oss then\qss $d\fff(Q)\qff =\qff 1$\qss 
for every component\dss $Q$\dss of\qss $S\cutc$\qss which is neither a disc with two holes\qss
({\fff}i.e.\qss a sphere with\dss $3$\dss holes),\oss
not a torus with\dss $1$\dss hole.\oss
But\qss $d\fff(Q)\qff =\qff 1$\qss if and only if\trs $Q$\trs is a sphere with\dss $4$\dss holes
or a torus with\dss $2$\dss holes.\oss
The theorem follows.\oss  \eproof

\mysection{Dehn\qss and\qss Dehn--Johnson\qss twists\qss 
in\qss Torelli\qss groups\halfff:\oss I}{simplest-twists-sufficient}

\myitpar{Commutants and bicommutants.} Let\qss $G$\qss be a group,\oss
and let\qss $X\dff\subset\dff G$\dnsp.\oss
The\qss \emph{commutant}\qss $X'$\qss of\dss $X$\dss is the set of all elements of\qss $G$\qss
commuting with all elements of\dss $X$\snsp.\oss
It is a subgroup of\qss $G$\dnsp.\oss 
The\qss \emph{bicommutant}\qss $X''$\qss of\dss $X$\dss is the commutant of\dss $X'$\snsp.\oss

The\qss \emph{commutant}\qss $g'$\qss of an element\qss $g\dff\in\dff G$\qss
is the set of all elements of\qss $G$\qss
commuting with\dss $g$\snsp.\oss
In other terms,\oss it is the commutant of the one-element subset\qss $\{\dff g\dff\}$\dnsp.\oss
The\qss \emph{bicommutant}\qss $g''$\qss of\dss $g$\dss is the commutant of\qss $g'$\snsp.\oss

In the rest of this section we will consider the commutants and bicommutants in\qss $\tors$\dnsp.\oss

\mypar{Theorem.}{forcing-multi-twist} \emph{Suppose that\oss 
$\genus(S)\qff\geq\qff 3$\snsp.\qff\oss
Suppose that\oss $f\dff\in\dff \tors$\dnsp.\off\oss 
If\oss
\vspace*{-\bigskipamount}
\begin{itemize}
\item $f$\oss belongs to an abelian subgroup of\oss 
$\tors$\oss of\qss rank\oss $2\dff\gen\qff -\qff 3$\nsp,\oss and\oss
\item $f''$\oss does not contain abelian subgroups\trs of\qss rank\qss $2$\nsp,\off\oss
\end{itemize}
\vspace*{-\bigskipamount} 
then\oss $f$\oss is a Dehn multi-twist\halfff.\oss}

\proof\qss Let\trs $\mathfrak{A}$\trs be an abelian subgroup of\qss 
$\tors$\qss containing\dss $f$\snsp.\oss
By Theorem\qss \ref{abelian-reduction}\qss there exists 
a pure reduction system\dss $c$\dss for\trs $\mathfrak{A}$\snsp.\oss
Then\dss $c$\dss is also a pure reduction system for\qss $f$\snsp,\oss 
and hence there is a diffeomorphism\qss $F\fff\colon S\ttoo S$\qss
representing\qss $f$\qss and such that\dss $c$\dss 
is a pure reduction system for\qss $F$\snsp.\oss
\vspace*{-\medskipamount}
\begin{quote}
\textsc{Claim.}\oss Let\qss $Q$\qss be a component of\qss $S\cutc$\qss
such that the isotopy class of\qss $F_Q$\qss is pseudo-Anosov,\qff\oss
and let\qss $C$\qss be a component of\pss $p\ccut c\dff(\partial Q)$\nsp.\qff\oss
Then the Dehn twist\qss $t_{\dff C}$\qss commutes with every element of the
commutant\qss $f'$\qss of\qss $f$\qss in\qss $\tors$\dnsp.\oss
\end{quote}
\vspace*{-\medskipamount}
\emph{Proof of the Claim.}\oss Let\dss $c_{\dff f}$\dss 
be a minimal pure reduction system for\qss $f$\qss contained in\dss $c$\snsp.\oss 
Suppose that\qss $g\dff\in\dff f'$\dnsp,\qff\oss
i.e.\qss that\qss $g\dff\in\dff \tors$\qss and\qss $g$\qss commutes with $f$\snsp.\qff\oss
Let\qss $G$\qss be a diffeomorphism of\qss $S$\qss representing\qss $g$\snsp.\qff\oss
Then\qss $G\fff(c_{\dff f})$\qss is a minimal pure reduction system for\oss 
$g\dff f\dff g^{\fff -\dff 1}\off =\off f$\nsp,\oss
and\dss hence\pss $G\fff(c_{\dff f})$\pss is\dss isotopic\dss to\qss $c_{\dff f}$\nsp.\oss
Replacing\qss $G$\qss by an isotopic diffeomorphism,\oss if necessary,\oss
we may assume that\pss  $G\fff(c_{\dff f})\off =\off c_{\dff f}$\dnsp.\oss 

Since the isotopy class of\qss $F_Q$\qss is pseudo-Anosov,\qff\oss
Lemma\qss \ref{boundary-pseudo-anosov}\qss implies that\pss 
$p\ccut c\dff(\partial Q)$\pss is contained in\dss $c_{\dff f}$\nsp.\oss
Hence Lemma\qss \ref{pure-abc}\qss implies that\qss 
$G$\qss leaves every component of\qss $p\ccut c\dff(\partial Q)$\qss invariant\halfff.\qff\oss
In particular\halfff,\oss $G\fff(C)\off =\off C$\dnsp,\oss 
and hence\oss
\[
\quad
g\qff t_{\dff C}\qff g^{\fff -\dff 1}\qff =\off t_{\dff C}\fff
\]
({\fff}because\qss $G$\qss is orientation-preserving).\oss
In other terms,\oss $t_{\dff C}$\qss commutes with\qss $g$\snsp.\oss  
Since\qss $g$\qss is an arbitrary element of\qss $f'$\snsp,\qff\oss
this proves the claim.\oss  \esubproof

Suppose that\qss $f$\qss is not a Dehn multi-twist\halfff.\qff\oss
Then there is a component\qss $Q$\qss of\qss $S\cutc$\qss such that
the isotopy class of\qss $F_Q$\qss is pseudo-Anosov.\qff\oss
Theorem\qss \ref{lower-rank-components}\qss implies that\qss $Q$\qss 
is either a sphere with\trs $3$\trs or\trs $4$\trs holes,\oss
or a torus with\trs $1$\trs or\trs $2$\qss holes.\oss

Since there are no pseudo-Anosov isotopy classes on a sphere with\dss $3$\dss holes,\oss
$Q$\qss cannot be a sphere with\dss $3$\dss holes.\oss
If\qss $Q$\qss is a torus with one hole,\oss
then every pseudo-Anosov isotopy class acts non-trivially on\qss $H_{\fff 1}(Q)$\dnsp.\oss
On the other hand,\oss in this case the inclusion homomorphism\qss 
$H_{\fff 1}(Q)\ttoo H_{\fff 1}(S)$\qss is injective.\oss
Since the isotopy class of\qss $F_Q$\qss is pseudo-Anosov,\oss
in this case\dss $f$\dss acts non-trivially on\qss $H_{\fff 1}(S)$\dnsp,\oss
in contradiction with\qss $f\dff\in\dff \tors$\dnsp.\oss
Therefore\qss $Q$\qss cannot be a torus with\dss $1$\dss hole either\halfff.\qff\oss
It follows that\qss $Q$\qss is either a sphere with\dss $4$\dss holes,\oss
or a torus with\dss $2$\dss holes.\oss

Suppose that\qss $Q$\qss is a torus with\dss $2$\dss holes.\oss
If\qss $p\cutc$\qss restricted to\qss $Q$\qss is not an embedding\halfff,\oss
then\qss $p\cutc$\qss maps both components of\qss $\partial Q$\qss onto 
the same circle in\qss $S$\dnsp.\oss
In this case the image\qss $p\cutc\fff(Q)$\qss is a 
closed subsurface of\qss $S$\qss and hence is equal to\qss $S$\dnsp.\oss
It follows that in this case\qss $S$\qss is a surface of genus\qss $2$\dnsp,\oss
contrary to the assumption.\oss
Therefore\qss $p\cutc$\qss embeds\qss $Q$\qss into\qss $S$\qss and we can
consider\qss $Q$\qss as a subsurface of\qss $S$\dnsp.\oss
Let\qss $C\fff,\pff D$\qss be the two boundary components of\qss $Q$\dnsp.\oss
Since the isotopy class of\qss $F_Q$\qss is pseudo-Anosov,\qss
the Dehn twists\qss $t_{\dff C}$\qss and\qss $t_{\dff D}$\qss commute with all 
elements of the commutant\qss $f'$\dnsp.\oss
While these Dehn twists themselves do not belong to\qss $\tors$\dnsp,\oss 
the product\oss
\[
\quad
t\off =\off t_{\dff C}\cdot\dff t_{\dff D}^{\dff -\halfff 1}
\]
is a Dehn-Johnson twist and hence belongs to\qss $\tors$\dnsp.\oss
Therefore both\qss $f$\qss and\qss $t$\qss belong to the bicommutant\qss $f''$\dnsp.\oss
Moreover\halfff,\oss by the classification of abelian subgroups of\qss $\mms$\qss
the elements\qss $f\fff,\pff t$\qss generate a free abelian group of rank\dss $2$\dnsp.\oss
Since\qss $f''$\qss cannot contain such a subgroup by the assumption,\dff\oss
it\dss follows\sss that\qss $Q$\qss is not a torus with\dss $2$\dss holes.\oss

The only possibility that remains is that\dss $Q$\dss is a sphere with\dss $4$\dss holes.\qff\oss
Theorem\qss \ref{lower-rank-necklaces}\qss implies that in this case either\qss 
$p\cutc$\qss restricted to\dss $Q$\dss is not an embedding\halfff,\oss
or\qss $\ccomp Q$\qss is not connected.\oss

If\qss $p\cutc$\qss restricted to\qss $Q$\qss is not an embedding\halfff,\oss
then\qss $p\cutc$\qss maps two components of\qss $\partial Q$\qss onto 
the same circle in\qss $S$\dnsp.\qff\oss
If\qss $p\cutc$\qss also maps the two other components of\qss $\partial Q$\qss
onto the same circle,\oss then the image\qss $p\cutc\fff(Q)$\qss is a 
closed subsurface of\qss $S$\qss and hence is equal to\qss $S$\dnsp.\oss
It follows that in this case\qss $S$\qss is a surface of genus\qss $2$\dnsp,\oss
contrary to the assumption.\oss
If\qss $p\cutc$\qss maps two other components of\qss $\partial Q$\qss to two
different circles\qss $C\fff,\pff D$\dnsp,\oss
then the image\qss $p\cutc\fff(Q)$\qss is a torus with\dss $2$\dss holes,\oss
and\qss $C\fff,\pff D$\qss are its boundary components.\oss
By arguing exactly as in the case of a torus with\dss $2$\dss holes,\oss
we conclude that in this case\qss $f''$\qss contains 
a free abelian group of rank\dss $2$\dss
having\qss $f$\qss and the Dehn-Johnson twist about 
the pair\qss $C\fff,\pff D$\qss as its free generators.\oss
This contradicts to the assumptions of the theorem,\oss
and hence\qss $p\cutc$\qss actually embeds\qss $Q$\qss into\qss $S$\dnsp.\oss

It remains to consider the case when\qss $Q$\qss 
is a sphere with\dss $4$\qss holes and\qss $Q$\qss is a subsurface of\qss $S$\dnsp.\qff\oss
As we already noted,\oss in this case\dss Theorem\qss \ref{lower-rank-necklaces}\qss
implies that\qss $\ccomp Q$\qss is not connected.\qff\oss
Therefore\qss $\ccomp Q$\qss consists of\pss 
$2$\dnsp,\oss $3$\dnsp,\oss or\qss $4$\pss components.\oss
In the first case the boundary of each of two components of\qss $\ccomp Q$\qss
is a bounding pair in\qss $S$\dnsp.\oss
As in the case of the torus with\dss $2$\dss holes,\oss
this implies that\qss $f''$\qss contains an abelian group of rank\dss $2$\dnsp.\oss
In the other two cases there is a component of\qss $\ccomp Q$\qss 
with only\dss $1$\dss boundary component.\qff\oss
Let\qss $C$\qss be this component\halfff.\qff\oss
Then\qss $C$\qss is a bounding circe,\oss and hence the Dehn twist\pss $t_{\dff C}$\pss
about\qss $C$\qss belongs to\qss $\tors$\dnsp.\oss
On the other hand,\oss $t_{\dff C}$\pss belongs to the bicommutant\pss $f''$\pss
by the above Claim.\qff\oss
Therefore,\oss in this case\qss $f''$\qss also contains 
a free abelian subgroup of rank\dss $2$\snsp,\oss
contrary to the assumption.\oss
 
It follows there is no component\qss $Q$\qss of\qss $S\cutc$\qss
such that the isotopy class of\qss $F_Q$\qss is pseudo-Anosov.\oss
Hence all diffeomorphisms\qss $F_Q$\qss are isotopic to the identity,\oss
and hence\qss $f$\qss is a product of Dehn twist about components of\dss $c$\dnsp.\oss  \eproof

\mypar{Theorem.}{bounding-twists} \emph{Under the assumption of
Theorem\qss \ref{forcing-multi-twist}\qss
$f$\qss is a non-zero power of either a Dehn twist about a separating circle,\oss
or a Dehn--Johnson twist about a bounding pair\halfff.\oss}

\proof\qss By Theorem\qss \ref{forcing-multi-twist},\oss 
$f$\qss is a Dehn multi-twist\halfff.\oss
Therefore,\oss $f$\qss has the form\qss (\ref{multi-twist-definition})\qss
for some one-dimensional submanifold\dss $c$\dss of\qss $S$\qss and some
integers\qss $m_{\dff O}$\snsp,\oss where\qss $O$\qss runs over the components of\dss $c$\snsp.\oss
Without any loss of generality we may assume that\dss $c$\dss is a system of circles
and that all integers\qss $m_{\dff O}\qff \neq\qff 0$\snsp.\oss
Then\dss $c$\dss is a minimal pure reduction system for\qss $f$\snsp.\oss

Let\qss $g\dff\in\dff \tors$\qss and let\qss $G$\qss 
be a diffeomorphism representing\qss $g$\snsp.\oss
If\qss $g\dff\in\dff f'$\snsp,\oss then\qss
$g\dff f\dff g^{\fff -\dff 1}\off =\off f$\nsp,\oss
and\pss $G\fff(c)$\pss is\dss isotopic\dss to\dss $c$\dss
because\dss $c$\dss is a minimal reduction system for\qss $f$\qss
(see Section\qss \ref{pure-reduction}).\oss
Replacing\qss $G$\qss by an isotopic diffeomorphism,\oss if necessary,\oss
we may assume that\pss  $G\fff(c)\off =\off c$\dnsp.\oss
Then by Lemma\qss \ref{pure-abc}\qss $G$\qss leaves every component of\dss $c$\dss invariant\halfff.\oss
It follows that\qss
$g\dff t_{\dff C}\dff g^{\fff -\dff 1}\off =\off t_{\dff C}$\qss
for every component\qss $C$\qss of\dss $c$\snsp.\oss
It follows that\qss $g$\qss commutes with all elements of\qss $\mathcal{T}\fff(c)$\dnsp.\oss
In particular\halfff,\oss $g$\qss commutes with all elements of\oss
$\mathcal{T}\fff(c)\fff\cap\dff\tors$\dnsp.\oss 
Since\qss $g$\qss is an arbitrary element of\qss $f'$\snsp,\oss
it follows that\qss
\[
\quad
\mathcal{T}\fff(c)\fff\cap\trf\tors\qff\subset\qff f''\fff.
\]
But Theorem\qss \ref{multi-twists-subgroups}\qss implies that\oss
$\rank \mathcal{T}\fff(c)\fff\cap\dff\tors\qff \geq\qff 2$\oss
unless\dss $c$\dss is either a separating circle,\oss
or the union of two circles forming a bounding pair\halfff.\oss
In the first case\qss $f$\qss is a power of the Dehn twist about this circe,\oss
and in the second case\qss $f$\qss is a power of a 
Dehn--Johnson twist about this bounding pair\halfff.\oss  \eproof

\mysection{Dehn\qss and\qss Dehn--Johnson\qss twists\qss 
in\qss Torelli\qss groups\halfff:\oss II}{simplest-twists-necessary}

\vspace*{\bigskipamount}
By Theorem\qss \ref{bounding-twists}\qss the two condition of
Theorem\qss \ref{forcing-multi-twist}\qss are sufficient for an element\dss $f$\dss of\qss $\tors$\dss
to be a non-zero power of either a Dehn twist about a separating circle,\oss
or a Dehn--Johnson twist about a bounding pair\halfff.\oss
This section is devoted to a proof that these conditions are necessary.\oss
See Theorem\qss \ref{simple-twists-satisfy}.\oss
The ideas of this proof are essentially the same as the ideas of the original proof
of B.\qss Farb and the author\halfff,\oss announced in\oss \cite{fi}.\oss
Naturally,\oss it was adapted to the context of the present paper
and differs from the original proof in details.

Theorems\qss \ref{bounding-twists}\qss and\qss \ref{simple-twists-satisfy}\qss
together provide an algebraic characterization of
non-zero powers of Dehn twists about separating circles and
Dehn--Johnson twists about a bounding pairs in terms of group structure of\qss $\tors$\dnsp.\oss
In order to distinguish between non-zero powers 
of Dehn twists about separating circles and
Dehn--Johnson twists about a bounding pairs,\oss 
one needs to use an additional algebraic condition.\oss 
See\oss \cite{fi},\oss Proposition\qss 9.\oss

\myitpar{A theorem of Thurston.} By a well known theorem of W. Thurston
one can construct pseudo-Anosov isotopy classes of diffeomorphisms of a surface $Q$
by taking the isotopy classes of various products of powers of twist diffeomorphism of\dss $Q$\dnsp.\oss
If one uses only twist diffeomorphisms about circles bounding in $Q$ 
subsurfaces with $1$ boundary component\halfff,\oss
then these products diffeomorphisms act trivially on\dss $H_{\fff 1}(Q)$\dnsp.\oss
This is how Thurston constructed the first examples of pseudo-Anosov isotopy classes
acting trivially on\dss $H_{\fff 1}(Q)$\dnsp,\oss
and,\oss in particular\halfff,\oss the first examples of pseudo-Anosov elements of Torelli groups\qss
(Thurston did not phrased his results in terms of\dss Torelli groups).\oss

In order for this construction to apply,\pss it is sufficient for $Q$ to contain circles
bounding in $Q$ subsurfaces with $1$ boundary component\halfff.\oss
This excludes only surfaces of genus $0$ and 
surfaces of genus $1$ with $1$ boundary component\halfff.\oss
When this construction applies,\oss it leads to many examples of pseudo-Anosov isotopy classes.\oss
In particular\halfff,\oss it leads to examples of non-commuting pseudo-Anosov isotopy classes.\oss

Suppose now that $Q$ is a subsurface of $S$\dnsp.\oss 
One can modify Thurston's examples by replacing the twist diffeomorphisms of $Q$ by
the twist diffeomorphisms of $S$ about the same circles.\oss
Then instead of a diffeomorphisms of $Q$ we will get diffeomorphisms of $S$ equal to
the identity on $\ccomp Q$ and such that the induced diffeomorphisms of $Q$ are pseudo-Anosov.\oss
If we use only twist diffeomorphisms about circles bounding in $Q$ 
subsurfaces with $1$ boundary component\halfff,\oss
then these diffeomorphisms of $S$ will act trivially on\dss $H_{\fff 1}(S)$\dnsp.\oss
In other terms,\oss their isotopy classes will belong to\dss $\tors$\dnsp.\oss 

It follows that if a subsurface $Q$ of $S$ is neither a surface of genus $0$ nor
a surface of genus $1$ with $1$ boundary component\halfff,\oss
then there are diffeomorphisms of $S$ equal to the identity on $\ccomp Q$\snsp,\oss
such that their isotopy classes belong to\dss $\tors$\dnsp,\oss 
and such that the induced diffeomorphisms of $Q$ 
belong to pseudo-Anosov isotopy classes.\oss
Moreover\halfff,\oss there are\qss (pairs of\dff)\qss 
such diffeomorphisms with non-commuting isotopy classes.\oss

\mypar{Lemma.}{torus} \emph{Let\dss $C$\dss be a circle on\qss $S$\qss
separating\qss $S$\qss into two parts\qss $Q\fff,\pff R$\qss 
having\dss $C$\dss as their common boundary.\oss
Suppose that\qss $G$\qss is a diffeomorphism of\qss $S$\qss 
leaving each of these parts invariant and such that its
isotopy class belongs to\qss $\tors$\dnsp.\oss
If\qss $R$\qss is a torus with\dss $1$\dss hole,\oss 
then the restriction\qss $G_{\fff R}$\qss is isotopic to the identity.\oss}

\proof\qss Since\qss $\partial R$\qss is a separating circle,\oss
the homology group\qss $H_{\fff 1}(R)$\qss is a direct summand of\qss $H_{\fff 1}(S)$\dnsp.\oss
By the assumption,\pss $G$\qss acts trivially on\qss $H_{\fff 1}(S)$\dnsp.\oss
It follows that\qss $G_{\fff R}$\qss acts trivially on\qss $H_{\fff 1}(R)$\dnsp.\oss
By the classification of diffeomorphisms of a torus with $1$ hole,\oss
this implies that\qss $G_{\fff R}$\qss is isotopic to the identity.\oss  \eproof

\mypar{Lemma.}{pa-parts} \emph{Let\dss $c$\dss be a system of circles on\qss $S$\qss
separating\qss $S$\qss into two parts having\dss $c$\dss as their common boundary.\oss
If both these parts have genus\dss $\geq\dff 1$\snsp,\oss 
then there is a  subset\pss $X\qff\subset\qff \tors$\pss with the commutant\qss $X'$\qss 
equal to\oss 
$\displaystyle 
\mathcal{T}\fff(c)\fff\cap\trf\tors$\dnsp.\oss}

\proof\qss Let\dss $Q\fff,\pff R$\dss be the two parts 
into which $c$\ divides\dss $S$\dnsp.\oss
Let us consider diffeomorphisms of\qss $S$\qss
equal to the identity on\dss $R$\dss
and such that the isotopy class of the induced diffeomorphism of\qss $Q$\qss 
is pseudo-Anosov.\oss
Suppose that\dss $Q$\dss is not a torus with $1$ hole.\oss
Then Thurston's construction leads,\pss in particular\halfff,\pss to two such diffeomorphisms\qss 
$F_{\fff 1}\fff,\pff F_{\fff 2}$\qss having the additional property 
that their isotopy classes\qss 
$f_{\fff 1}\fff,\pff f_{\fff 2}$\qss belong to\qss $\tors$\qss and do not commute.\oss
If\qss $R$\qss is also not a torus with $1$ hole,\oss
then there are also diffeomorphisms\qss 
$G_{\fff 1}\fff,\pff G_{\fff 2}$\qss having the same properties,\oss
but with the roles of\dss $Q$\dss and\dss $R$\dss interchanged.\oss
Let\dss $g_{\fff 1}\fff,\pff g_{\fff 2}$\dss be their isotopy classes,\oss
and let\oss 
\[
\quad
X\off =\off \{\dff f_{\fff 1}\fff,\pff f_{\fff 2}\fff,\pff g_{\fff 1}\fff,\pff g_{\fff 2}\dff\}\dff.
\]
Then Dehn twists about components of\dss $c$\dss 
commute with all elements of\dss $X$\qss
and\dss hence\oss
$\displaystyle
\mathcal{T}\fff(c)\fff\cap\trf\tors\qff\subset\qff X'$\dnsp.\oss
Let us prove the opposite inclusion\oss
$\displaystyle
X'\qff\subset\qff \mathcal{T}\fff(c)\fff\cap\trf\tors$\dnsp.\oss

Let\qss $f\dff\in\dff X'$\dnsp.\oss
Then the subgroup\qss $\mathfrak{A}$\qss of\qss $\tors$\qss generated by\dss $f$\dss and,\oss 
say,\oss $f_{\fff 1}\fff,\pff g_{\fff 1}$\qss is abelian.\oss
By Theorem\qss \ref{abelian-reduction}\qss there exists a pure reduction
system\dss $c_{\fff 0}$\dss for\qss $\mathfrak{A}$\snsp.\oss
Then\dss $c_{\fff 0}$\dss is also a pure reduction system for\dss $f_{\fff 1}$\snsp.\oss
By Lemma\qss \ref{boundary-pseudo-anosov}\qss $c$\dss is contained up to isotopy 
in any pure reduction system for\dss $f_{\fff 1}$\snsp,\oss
and hence up to isotopy\dss $c$\dss is contained in\dss $c_{\fff 0}$\snsp.\oss
Therefore we may assume that\qss $c_{\fff 0}\qff \supset\qff c$\snsp.\oss
Then\qss $c_{\fff 0}\qff =\qff c$\qss because diffeomorphisms\qss 
$F_{\fff 1}\fff,\pff G_{\fff 1}$\qss cannot leave invariant 
any system of circles in\qss $Q\fff,\pff R$\qss respectively.\oss 

Now the description of abelian subgroups of\qss $\tors$\sss 
from Section\qss \ref{pure-reduction}\qss implies that\qss
$\mathfrak{A}$\qss is contained in the group
generated by\qss $\mathcal{T}\fff(c)\fff\cap\trf\tors$\qss 
and\qss $f_{\fff 1}\fff,\pff g_{\fff 1}$\snsp.\oss
It follows that\oss 
\begin{equation}
\label{f-general}
\quad
f\off =\off t\cdot f_{\fff 1}^{\trf m}\cdot g_{\fff 1}^{n}\dff
\end{equation}
for some\qss $t\qff\in\qff \mathcal{T}\fff(c)\fff\cap\trf\tors$\qss
and\qss $m\fff,\pff n\dff\in\dff \zzz$\snsp.\oss
Since\qss $f_{\fff 2}$\qss commutes with\qss $t$\nsp,\qss $g_{\fff 1}$\nsp,\pss and\pss $f$\nsp,\oss 
(\ref{f-general})\qss implies that\qss $f_{\fff 2}$\qss commutes with\qss $f_{\fff 1}^{\trf m}$\nsp.\oss 
If\qss $m\qff \neq\qff 0$\snsp,\oss this implies that\qss $f_{\fff 2}$\qss 
commutes with\qss $f_{\fff 1}$\qss contrary to the assumption.\oss
Therefore\qss $m\qff =\qff 0$\snsp.\oss
By a completely similar argument\qss $n\qff =\qff 0$\qss and hence\oss
\[
\quad
f\off =\off t\qff\in\qff \mathcal{T}\fff(c)\fff\cap\trf\tors\dff.
\]
This proves the lemma in the case when neither\dss $Q$\dnsp,\pss nor 
$R$\dss is a torus with $1$ hole.\oss 

Suppose now that\qss one of the surfaces\qss $Q\fff,\pff R$\qss is a torus
with $1$ hole,\oss but the other is not\halfff.\oss
We may assume that\qss $R$\qss is a torus with $1$ hole 
and\dss $Q$\dss is not\halfff.\oss
Then Thurston's construction applies to\dss $Q$\dss and leads to diffeomorphisms\qss 
$F_{\fff 1}\fff,\pff F_{\fff 2}$\qss with the same properties as above.\oss
Let\qss $f_{\fff 1}\fff,\pff f_{\fff 2}$\qss be their isotopy classes,\oss
and let\oss 
\[
\quad 
X\off =\off \{\dff f_{\fff 1}\fff,\pff f_{\fff 2}\dff\}\dff.
\]
Then the Dehn twist about the circle\qss $\partial Q\qff =\qff \partial R\qff =\qff c$\qss
belongs to\dss $X'$\dss and generates\oss
$\displaystyle
\mathcal{T}\fff(c)
\off =\off
\mathcal{T}\fff(c)\fff\cap\trf\tors$\dnsp.\oss
Therefore\oss
$\displaystyle 
\mathcal{T}\fff(c)\fff\cap\trf\tors\qff\subset\qff X'$\dnsp.\oss

In order to prove the opposite inclusion,\oss
consider an arbitrary\qss $f\dff\in\dff X'$\qss and the group\qss
$\mathfrak{A}$\qss generated by\dss $f$\dss and\dss $f_{\fff 1}$\snsp.\oss
Arguing as above,\oss we see that there is a pure reduction system\dss $c_{\fff 0}$\dss 
for\qss $\mathfrak{A}$\qss containing\dss $c$\snsp.\oss
In view of Lemma\qss \ref{torus},\pss 
this implies that every element\qss $\mathfrak{A}$\qss can
be represented by a diffeomorphism\dss $G$\dss leaving\dss $R$\dss invariant and such that
the restriction\qss $G_{\fff R}$\qss is isotopic to the identity.\oss
It follows that\qss $\mathfrak{A}$\qss is contained in the group
generated by\qss $\mathcal{T}\fff(c)$\qss 
and\qss $f_{\fff 1}$\qss and hence\oss 
\begin{equation}
\label{f-torus}
\quad
f\off =\off t\cdot f_{\fff 1}^{\trf m}
\end{equation}
for some\qss $t\qff\in\qff \mathcal{T}\fff(c)$\qss
and\qss $m\dff\in\dff \zzz$\snsp.\qff\oss
Arguing as above,\oss we see that\qss $m\qff =\qff 0$\snsp.\oss
Hence\oss
\[
\quad
f\off =\off t\qff\in\qff \mathcal{T}\fff(c)\off =\off \mathcal{T}\fff(c)\fff\cap\trf\tors\dff.
\]
This proves the lemma in the case when one of the parts\dss $Q\fff,\pff R$\dss  
is a torus with $1$ hole.\oss

If each part\qss $Q\fff,\pff R$\qss is a torus with $1$ hole,\oss
then one can take\qss $X\qff =\qff \{\trf t_{\dff C}\dff\}$\snsp,\oss
where\qss $C\qff =\qff \partial Q\qff =\qff \partial R$\snsp.\oss
We leave the details of this case to the reader\halfff.\oss  \eproof

\mypar{Theorem.}{simple-twists-satisfy} \emph{If\qss $f$\qss is a non-zero power of either 
Dehn twist about a separating circle,\oss
or Dehn--Johnson twist about a bounding pair,\pss
then the\dss $f$\dss satisfies the conditions of 
Theorem\dss \ref{forcing-multi-twist}.\oss}

\proof\qss Suppose that\dss $f$\dss is a non-zero power of a Dehn twist about 
a separating circle\qss $c$\snsp.\qff\oss
Let\pss $Q_{\fff 1}$\dnsp,\qss $Q_{\fff 2}$\pss be the parts 
into which\dss $C$\dss divides\dss $S$\dnsp,\oss
and let\qss $\gen_{\fff 1}\qff =\qff \genus(Q_{\fff 1})$\dnsp,\oss
$\gen_{\fff 2}\qff =\qff \genus(Q_{\fff 2})$\dnsp.\oss
Then\qss $\gen_{\fff 1}\qff +\qff \gen_{\fff 2}\qff =\qff \gen$\dnsp,\oss
and the pair\qss $(S\fff,\pff C)$\qss is determined up to a diffeomorphisms by 
the numbers\qss $\gen_{\fff 1}\fff,\qss \gen_{\fff 2}$\nsp.\oss 

The system of circles\dss $s$\dss illustrated on Fig.\qss \ref{separating-circles}\dss
contains a circle dividing\qss $S$\qss into two parts of 
genus\qss $\gen_{\fff 1}\fff,\pff \gen_{\fff 2}$\qss for every pair\qss
$\gen_{\fff 1}\fff,\pff \gen_{\fff 2}$\qss such that\qss 
$\gen_{\fff 1}\qff +\qff \gen_{\fff 2}\qff =\qff \gen$\dnsp.\oss
We may assume that\qss $c$\qss is one of these circles.\oss
Then\qss $f\dff\in\dff \mathcal{T}\fff(s)\fff\cap\dff\tors$\dnsp.\oss
The system of circles\dss $s$\dss consists of\qss $2\dff\gen\qss -\qss 3$\qss components
and has no necklaces because all components of\dss $s$\dss are separating\halfff.\oss 
Therefore,\oss
$\rank\qff \mathcal{T}\fff(s)\fff\cap\dff\tors\qff =\qff 2\dff\gen\qff -\qff 3$\oss
by Theorem\qss \ref{multi-twists-subgroups}.\oss
It follows that\dss $f$\dss satisfies the first condition of 
Theorem\qss \ref{forcing-multi-twist}.\oss

By Lemma\qss \ref{pa-parts},\oss there is a  subset\pss 
$X\qff\subset\qff \tors$\pss with\oss 
$\displaystyle 
X'\off =\off \mathcal{T}\fff(c)$\dnsp.\oss
Since\qss $f\dff\in\dff \mathcal{T}\fff(s)\fff\cap\dff\tors$\dnsp,\oss
the element\dss $f$\dss commutes with all elements of\qss $X$\snsp.\oss
Therefore,\oss $X\qff\subset\qff f'$\qss
and hence\qss 
$\displaystyle
f''\qff\subset\qff X'\off =\off \mathcal{T}\fff(c)\fff\cap\trf\tors\off =\off \mathcal{T}\fff(c)$\dnsp.\oss
Since in this case\qss $\mathcal{T}\fff(c)$\qss is an infinite cyclic group,\oss
$f''$\qss does not contains free abelian groups of rank $2$\snsp.\oss
It follows that\qss $f$\qss satisfies the second condition of 
Theorem\qss \ref{forcing-multi-twist}.\oss

Suppose now that\qss $f$\qss is a non-zero power of a Dehn--Johnson twist about 
a bounding pair\qss $C\fff,\pff D$\snsp.\off\oss
Let\oss
$\displaystyle
c\off =\off C\dff\cup\dff D$\dnsp.\qff\oss
Let\pss $Q_{\fff 1}$\dnsp,\qss $Q_{\fff 2}$\pss be the parts 
into which\dss $c$\dss divides\dss $S$\dnsp,\oss
and let\qss $\gen_{\fff 1}\qff =\qff \genus(Q_{\fff 1})$\dnsp,\oss
$\gen_{\fff 2}\qff =\qff \genus(Q_{\fff 2})$\dnsp.\oss
Then\qss $\gen_{\fff 1}\qff +\qff \gen_{\fff 2}\qff =\qff \gen\qff -\qff 1$\snsp,\oss
and the pair\qss $(S\fff,\pff c)$\qss is determined up to a diffeomorphisms by 
the numbers\qss $\gen_{\fff 1}\fff,\qss \gen_{\fff 2}$\nsp.\oss 

The system of circles\dss $s'$\dss illustrated on Fig.\qss \ref{sep-non-sep-cicrles}\dss
contains a two circles forming a bounding pair and dividing\qss $S$\qss into two parts of 
genus\qss $\gen_{\fff 1}\fff,\pff \gen_{\fff 2}$\qss for every pair\qss
$\gen_{\fff 1}\fff,\pff \gen_{\fff 2}$\qss such that\qss 
$\gen_{\fff 1}\qff +\qff \gen_{\fff 2}\qff =\qff \gen\qff -\qff 1$\snsp.\oss
As above,\oss we may assume that\qss $C\fff,\pff D$\qss is one of these bounding pairs.\oss
Then\qss $f\dff\in \mathcal{T}\fff(s')\fff\cap\dff\tors$\dnsp.\oss
The system of circles\dss $s'$\dss consists of\qss $2\dff\gen\qss -\qss 2$\qss components
and\dss $1$\dss necklace.\oss 
Therefore,\oss
$\rank\qff \mathcal{T}\fff(s')\fff\cap\dff\tors\qff =\qff 2\dff\gen\qff -\qff 3$\oss
by Theorem\qss \ref{multi-twists-subgroups}.\oss
This proves that\qss $f$\qss satisfies the first condition
of Theorem\qss \ref{forcing-multi-twist}.\oss

By using Lemma\qss \ref{pa-parts}\qss in the same way as before,\oss
we see that\qss 
$\displaystyle
f''\qff\subset\qff X'\off =\off \mathcal{T}\fff(c)\fff\cap\trf\tors$\dnsp.\oss
In this case\qss $\mathcal{T}\fff(c)\fff\cap\trf\tors$\qss is an infinite
cyclic group generated by a Dehn--Johnson twist about the bounding pair\qss $C\fff,\pff D$\dnsp.\oss
It follows that\qss $f$\qss satisfies the second condition of 
Theorem\qss \ref{forcing-multi-twist}.\oss
This completes the proof.\oss  \eproof

\myappend{BP-necklaces}{necklaces}

\myitpar{BP-necklaces.} A one-dimensional 
closed submanifold $c$\qss of\qss $S$\qss 
is called a\pss \emph{BP-necklace}\pss 
if the cut surface\qss $S\cutc$\qss 
has at least\dss $2$\dss components and every component of\qss $S\cutc$\qss 
has exactly\dss $2$\dss boundary components.\oss

\renewcommand{\topfraction}{1.0}
\renewcommand{\textfraction}{0.0}

\begin{figure}[t]
\includegraphics[width=0.96\textwidth]{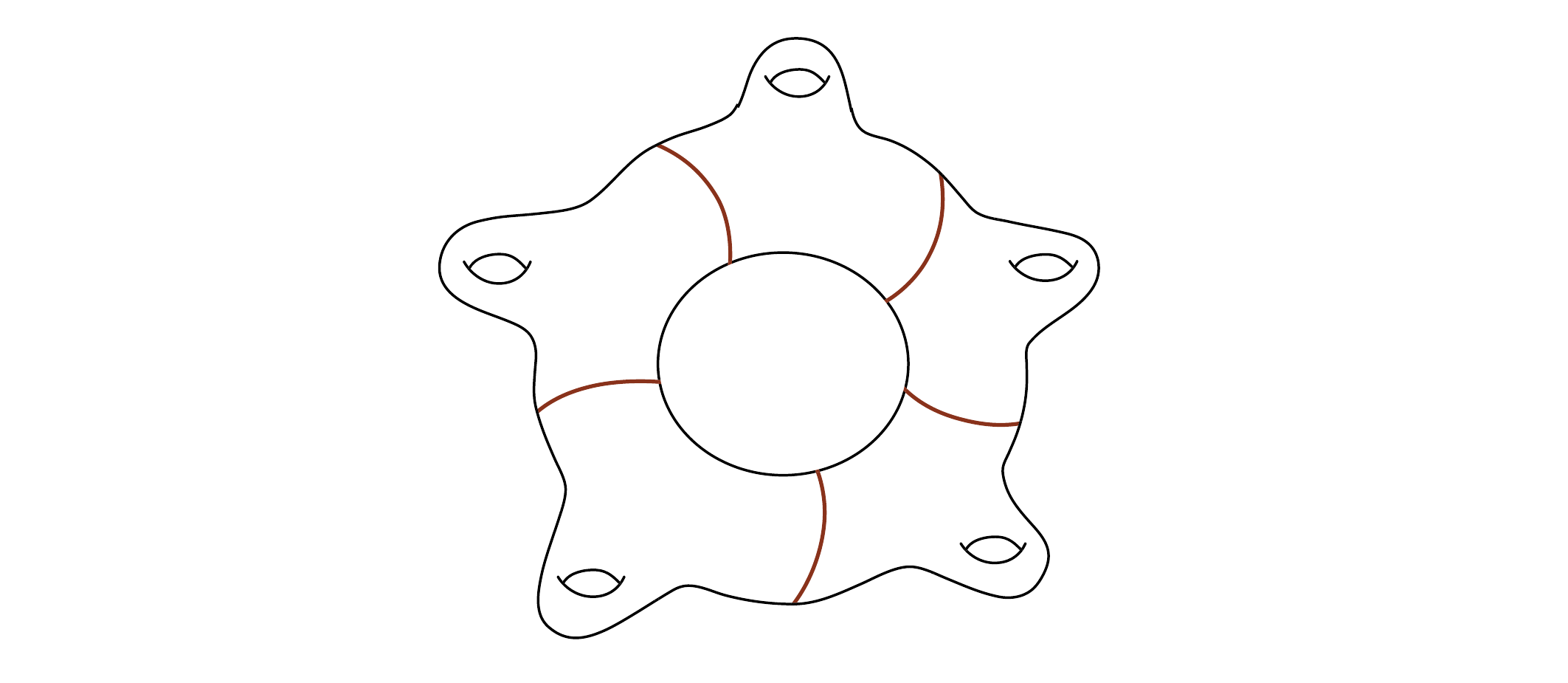}
\caption{A\dss BP-necklace.}
\label{bp-necklace}
\end{figure}

\begin{figure}[h]
\includegraphics[width=0.96\textwidth]{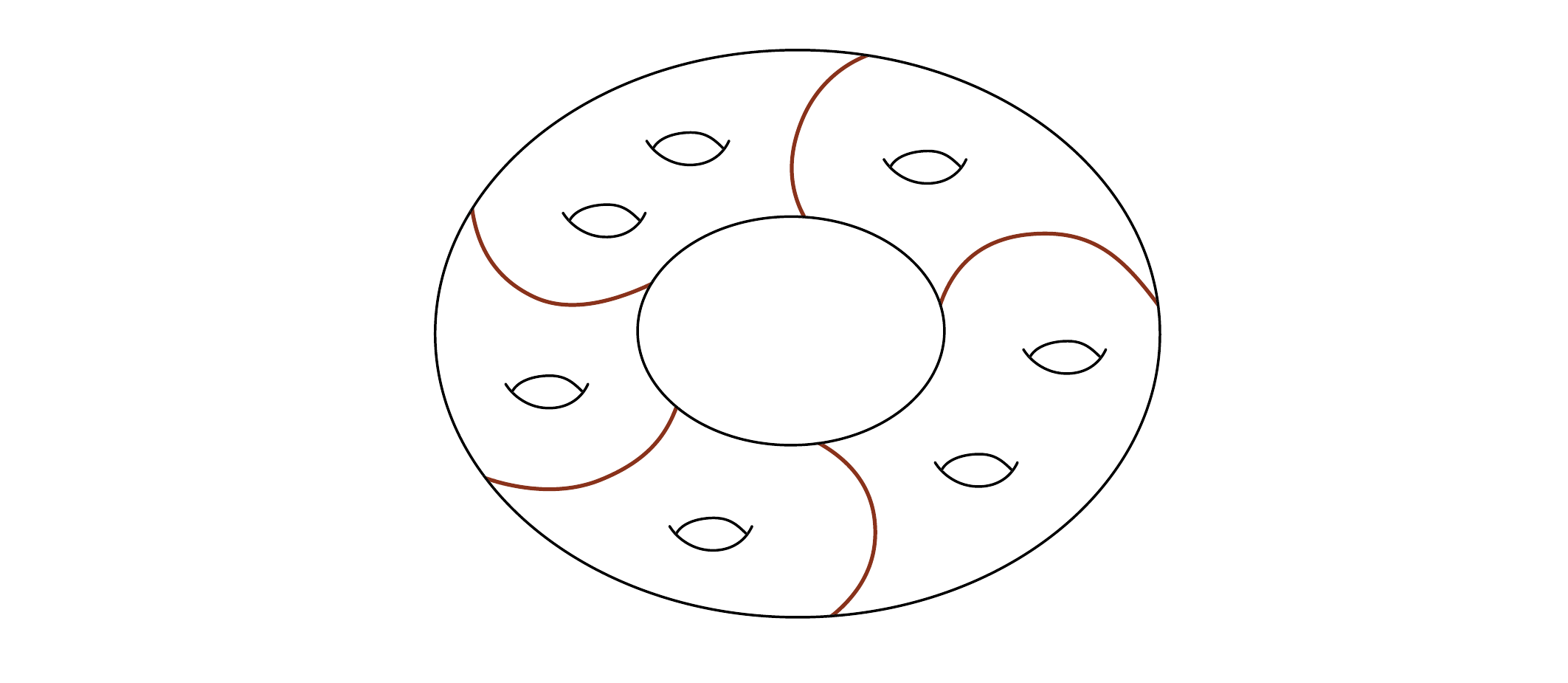}
\caption{Another\dss BP-necklace.}
\label{bp-necklace-general}
\end{figure}

An important example\sss $c$\sss of a BP-necklace on\dss $S$\dss is illustrated by Fig.\qss \ref{bp-necklace}.\oss
It is a union of\qss $\gen\qff -\qff 1$\qss circles.\oss
Every component of\qss $S\cutc$\qss is a torus with\dss $2$\dss holes.\oss
By adding to\dss $c$\dss for each component\qss $Q$\qss of\qss $S\cutc$\qss
a circle bounding a torus with\sss $1$\sss hole in\qss $Q$\dnsp,\oss
we get a system of circles\dss $s$\dss in\qss $S$\qss 
consisting of\qss $2\dff\gen\qff -\qff 2$\qss components and having exactly\dss 1\dss
necklace,\oss namely the set of components of\dss $c$\snsp.\oss
See Fig.\qss \ref{sep-non-sep-cicrles}.\oss
As we will see now,\oss every BP-necklace looks like the one on Fig.\dss \ref{bp-necklace},\oss
or\halfff,\oss more generally,\oss as the one on Fig.\dss \ref{bp-necklace-general}.\oss

\renewcommand{\topfraction}{1.0}
\renewcommand{\textfraction}{0.0}

\myitpar{A description of BP-necklaces.}
Let\dss $c$\dss be a BP-necklace.\oss
If\qss $p\cutc$\qss is not an embedding on a component\qss $Q$\qss of\qss $S\cutc$\snsp,\oss
then\qss $p\cutc$\qss maps both boundary circles of\qss $Q$\qss 
onto the same circle in\qss $S$\dnsp.\oss
In this case\qss $p\cutc\dff(Q)$\qss is a closed surface and hence is equal to\qss $S$\dnsp.\oss
Therefore\qss $Q$\qss is the only component of\qss $S\cutc$\snsp,\oss
contrary to the assumption.\oss
Hence\qss $p\cutc$\qss embeds all components of\qss $S\cutc$\qss and 
we may consider them as subsurfaces of\qss $S$\dnsp.\oss 

Every component of\dss $c$\dss is a component of the boundary of exactly two parts of\qss $S$\dnsp.\oss
Let\qss $Q_{\fff 0}$\qss be one of the parts of\dss $S$\dss with respect to\dss $c$\nsp,\oss
and let\qss $C_{\fff 0}$\qss be one of components of the boundary\dss $\partial Q_{\fff 0}$\snsp.\qff\oss
Let\qss $C_{\fff 1}$\qss be the other component of the boundary\qss $\partial Q_{\fff 0}$\dnsp,\oss
and\dss let\qss $Q_{\fff 1}$\qss be the second part of\qss $S$\qss
having\qss $C_{\fff 1}$\qss as a component of its boundary.\oss 
By continuing in this way,\oss one can consecutively number the components of\dss $c$\dss 
and the parts of\dss $S$\dss as
\[
\quad
C_{\fff 0}\fff,\off Q_{\fff 0}\fff,\quad 
C_{\fff 1}\fff,\off Q_{\fff 1}\fff,\quad
C_{\fff 2}\fff,\off Q_{\fff 2}\fff,\quad 
\ldots\fff,\off\off 
C_{\fff n}\fff,\off Q_{\fff n}\fff
\]
in such a way that\oss 
$\partial Q_{\fff i}\off =\off C_{\fff i}\cup C_{\fff i\qff +\qff 1}$\oss
if\oss $0\dff\leq i\dff\leq\dff n\qff -\qff 1$\dnsp,\oss
and\oss
$\partial Q_{\fff n}\off =\off C_{\fff n}\cup C_{\fff 0}$\nsp.\oss
In particular\halfff,\oss we see that all components of\dss $c$\dss
are homology equivalent\halfff.\oss

\myappar{Theorem.}{he-is-necklace} \emph{If a closed one-dimensional 
submanifold\dss $c$\dss of\qss $S$\qss has\qss $\geq\dff 2$\qss components
and all components of\dss $c$\dss are non-separating and homology equivalent\halfff,\oss
then\dss $c$\dss is a BP-necklace.\oss}

\proof If two disjoint non-separating circles\qss 
$C\fff,\pff C'$\qss are homology equivalent\halfff,\oss
then\qss $C\fff,\pff C'$\qss is a bounding pair\halfff.\oss
Indeed,\oss if the union $C\cup C'$ does not separate\dss $S$\dnsp,\oss 
then the classification of surfaces implies that 
the homology classes\qss $\hclass{C}$\snsp,\oss $\hclass{C'}$\qss
are linearly independent and hence\qss $C\fff,\pff C'$\qss 
cannot be homology equivalent\halfff.\oss
This observation implies the theorem in the case when\dss $c$\dss consists of\dss $2$\dss components.\oss
In order to prove the theorem in the general case,\oss
we will use the induction by the number of components of\dss $c$\dnsp.\oss  
Suppose that the theorem is already proved for submanifolds with\qss 
$\leq\dff n\qff -\qff 1$\qss components,\oss
and that\dss $c$\dss has\qss $n$\qss components.\oss  
Let\qss $C$\qss be a component of\dss $c$\dnsp,\oss
and let\qss $c_{\dff 0}\qff =\qff c\smallsetminus C$\dnsp.\oss 
By the inductive assumption,\oss $c_{\dff 0}$\qss is a BP-necklace.\oss
 
Since the circle\dss $C$\dss is disjoint from\dss $c_{\dff 0}$\nsp,\oss 
it is contained in a component\qss $Q$\qss of\qss $S\ccut c_{\dff 0}$\nsp.\oss 
If\qss $C$\qss is a non-separating in\qss $Q$\dnsp,\oss 
then\qss $C$\qss is not homology equivalent to components of\qss $\partial Q$\dnsp,\oss 
contrary to the assumption.\oss
Hence\qss $C$\qss separates\qss $Q$\qss into two parts.\oss 
If both components of\qss $\partial Q$\qss are contained in the same part,\oss 
then the other part is a subsurface of\qss $Q$\qss with boundary equal to\qss $C$\dnsp.\oss
This subsurface is also a subsurface of\qss $S$\dnsp.\oss 
Therefore in this case\qss $C$\qss bounds a subsurface of\qss $S$\dnsp,\oss
contrary to the assumption that\qss $C$\qss is non-separating.\oss
Therefore the two components of\qss $\partial Q$\qss 
are contained in different parts of\qss $Q$\dnsp.\oss
Hence\qss $C$\qss divides\qss $Q$\qss into two subsurfaces,\oss 
each of which has two boundary components.\oss 
One of these two boundary components is\qss $C$\dnsp,\oss  
and the other is a component of\qss $\partial Q$\dnsp.\oss 
It follows that\dss $c$\dss divides\qss $S$\qss into
components of\qss $S\ccut c_{\dff 0}$\qss different from\qss $Q$\qss and  
the two parts into which\qss $C$\qss divides\qss $Q$\dnsp.\oss  
Therefore\dss $c$\dss is a BP-necklace.\oss  \eproof

\myappar{Corollary.}{BP-and-necklaces} \emph{If\dss $s$\dss is a closed one-dimensional
submanifold of\qss $S$\dnsp,\oss then the union of circles in any necklace of\dss $s$\dss
is a BP-necklace.\oss}  \eproof  

%\vspace*{\bigskipamount}
%\asterism

%\newpage

\vspace*{\bigskipamount}

\begin{flushright}

June 24\fff,\oss 2016
 
https\halfff:/\!/\!nikolaivivanov.com

%E-mail\halfff:\oss radiantsadness@gmail.com

\end{flushright}

\end{document}